\documentclass[a4paper,12pt,reqno]{amsart}
\usepackage{amsmath,graphicx,graphics,amssymb,hyperref}
\newtheorem{theo}{Theorem}[section]
\newtheorem{lem}[theo]{Lemma}

\newcommand{\emailhref}[1]{\email{\href{#1}{#1}}}

\numberwithin{equation}{section}

\newcommand{\M}{\operatorname{M}}

\newcommand{\h}{\operatorname{H}}

\newcommand{\Z}{\mathbb{Z}}

\newmuskip\pFqskip
\mathchardef\pFcomma=\mathcode`, 

\begin{document}

\title{A short combinatorial proof of Di Francesco's conjecture on Aztec triangles}

\author[Seok Hyun Byun]{Seok Hyun Byun}\emailhref{sbyun@clemson.edu}
\address{School of Mathematical and Statistical Sciences, Clemson University, Clemson, South Carolina 29631, U.S.A.}
\thanks{S.H.B. was supported in part by an AMS-Simons Travel Grant.}

\author[Mihai Ciucu]{Mihai Ciucu}\emailhref{mciucu@iu.edu}
\address{Department of Mathematics, Indiana University, Bloomington, Indiana 47405, U.S.A.}
\thanks{M.C. was supported in part by Simons Foundation Collaboration Grant 710477.}



\begin{abstract} Di Francesco conjectured in 2021 that the number of domino tilings of a certain family of regions --- called Aztec triangles --- on the square lattice is given by a product formula reminiscent of the one giving the number of alternating sign matrices. This turned out to be a real challenge to prove without the use of computers --- each of the two existing proofs (one due to Koutschan, Krattenthaler and Schlosser, the other to Corteel, Huang and Krattenthaler) relies on substantial computer calculations which would be hard to check directly. In this paper we present a short combinatorial proof that relies on the second author's factorization theorem and complementation theorem for perfect matchings.

\end{abstract}

\maketitle

%
%

\section{Introduction}

\begin{figure}[t]
\centerline{
\hfill
{\includegraphics[width=0.20\textwidth]{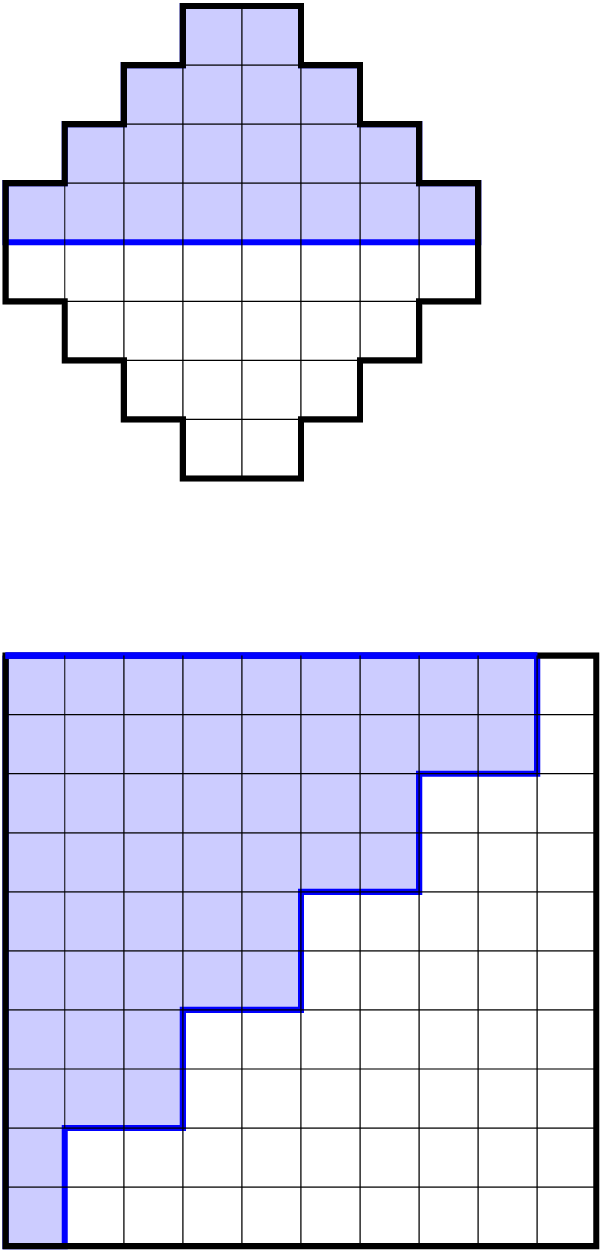}}
\hfill
{\includegraphics[width=0.18\textwidth]{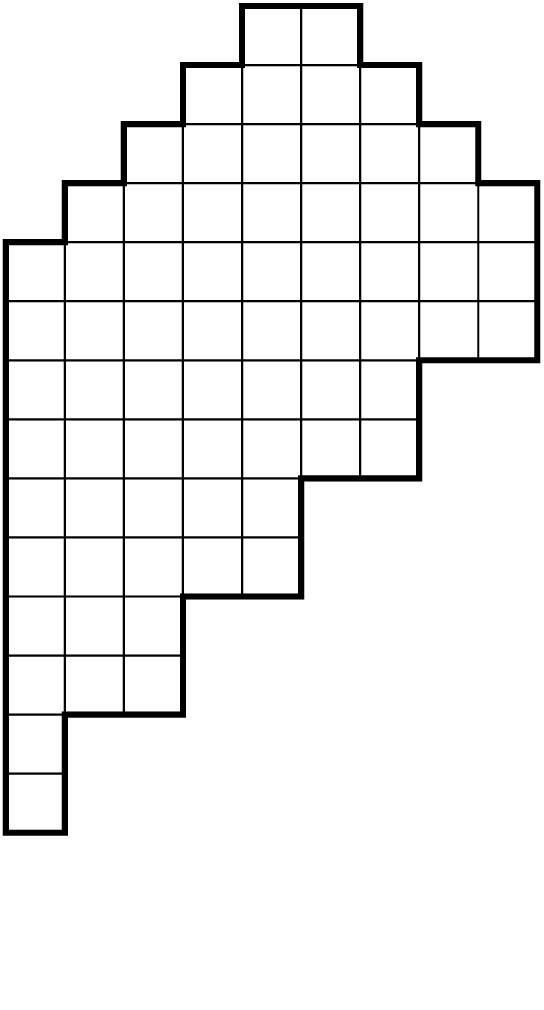}}
\hfill
{\includegraphics[width=0.18\textwidth]{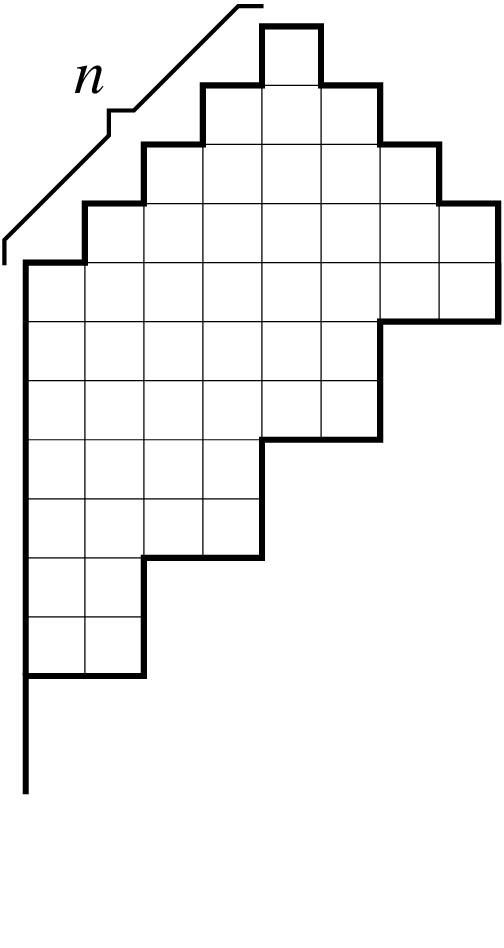}}
\hfill
}
\vskip0.1in
\caption{Half of the square of side length $2n$, and half of the Aztec diamond of order $n-1$, for $n=5$ (left). The region ${\mathcal T}_n$ for $n=5$ (center), and its planar dual graph (right).}
\vskip-0.1in
\label{faa}
\end{figure}

Let $S_{2n}$ be the lattice square of side length $2n$. Cut $S_{2n}$ into two congruent parts along its diagonal parallel to the first bisector by a zigzag line with steps of length two, leaving the bottom left unit square above the cut (this is illustrated on the bottom left in Figure \ref{faa}); denote the region above the cut by $HS_{2n}$. Let $AD_{n-1}$ be the Aztec diamond region of order $n-1$ ($AD_4$ is shown on the top left in Figure \ref{faa}), and denote by $HD_{n-1}$ its top half. Then the Aztec triangle of order $n$, denoted ${\mathcal T}_n$, is obtained by placing $HD_{n-1}$ on top of $HS_{2n}$ so that they are right-justified (see Figure \ref{faa}).

Di Francesco conjectured \cite[Conjecture 8.1]{DF2021} (see also \cite{DFG2020},  where Aztec triangles were first considered, and connections to the twenty-vertex model of statistical physics are explained) that the number of domino tilings of ${\mathcal T}_n$ is given by\footnote{ For a region $R$ on the square lattice, we denote by $\M(R)$ the number of domino tilings of $R$.}
\begin{equation}
\M({\mathcal T}_n)=
2^{n(n-1)/2}\prod_{i=0}^{n-1}\frac{(4i+2)!}{(n+2i+1)!},
\label{eaa}  
\end{equation}
a formula that, as he points out, is reminiscent of the one giving the number $A_n$ of alternating sign matrices of order $n$, which is $A_n=\prod_{i=0}^{n-1}\frac{(3i+1)!}{(n+i)!}$ (see \cite{Zeil,Kup,Fischer} for three different proofs).

Proving this formula without the use of computers turned out to be quite a challenge.
In \cite{Ciu2022} the second author made partial progress by showing that ${\mathcal T}_n$ is a divisor of a certain product similar to the right hand side of \eqref{eaa} (which gives the number of domino tilings of the cruciform regions introduced in that paper). The first proof of formula \eqref{eaa} was obtained by Koutschan, Krattenthaler and Schlosser~\cite{KKS2025}, an article which extends Koutschan's original proof \cite{Kou21} of formula \eqref{eaa}  (which boiled down to the evaluation of a certain determinant) to proving more general determinant identities. This work is based on the holonomic Ansatz developed by Zeilberger \cite{ZeilHol}, a computer-algebra-based approach that can be used to prove the evaluation of certain symbolic determinants in which the dimension of the square matrix is given by a symbolic parameter. As the authors state, due to their size, the results of the intermediate calculations cannot be displayed in their article (and are shown instead in an accompanying electronic material \cite{KKSelectr}).

The other existing proof of formula \eqref{eaa} is due to Corteel, Huang and Krattenthaler \cite{CHK2023}. This uses the more classical ``factor exhaustion'' method, which proves that a determinant is equal to a polynomial that factors into linear factors by providing independent linear combinations of the rows or columns that vanish for certain required specializations of the parameters. Much of this can be carried out directly, but for one essential family of linear factors the authors had to rely on substantial computer calculations, using Koutschan's {\tt HolonomicFunctions} {\it Mathematica} package \cite{KouHol} (see the Appendix in \cite{CHK2023}).

The purpose of this paper is to give a simple combinatorial proof of Di Francesco's formula~\eqref{eaa}, based on the second author's factorization theorem \cite{Ciu1997} and complementation theorem~\cite{Ciu1998}. It uses in an essential way the cruciform regions the second author considered in \cite{Ciu2022} and the formula found there giving the number of their domino tilings. The cruciform regions were arrived at by symmetrizing Di Francesco's Aztec triangles three times in succession. In essence, the proof presented here undoes these three steps in such a way that, due to a useful ``cancellation trick,'' at each step the number of tilings of the smaller regions can be expressed in terms of the number of tilings of the larger ones. The first of these steps was actually achieved by the second author in \cite{Ciu2022}; bringing this to completion fulfills the hope stated by one of the referees of that paper. 

\section{The proof}

We deduce formula \eqref{eaa} from two results on the number of perfect matchings of cruciform-shaped graphs. The first of these is the second author's \cite[Theorem 2.1]{Ciu2022}, which we recall here for convenience. Throughout this paper, instead of domino tilings of regions on the square grid, we will use the equivalent language of perfect matchings (for simplicity, often called just matchings) of their planar duals.

\begin{figure}[t]
\centerline{
{\includegraphics[width=0.43\textwidth]{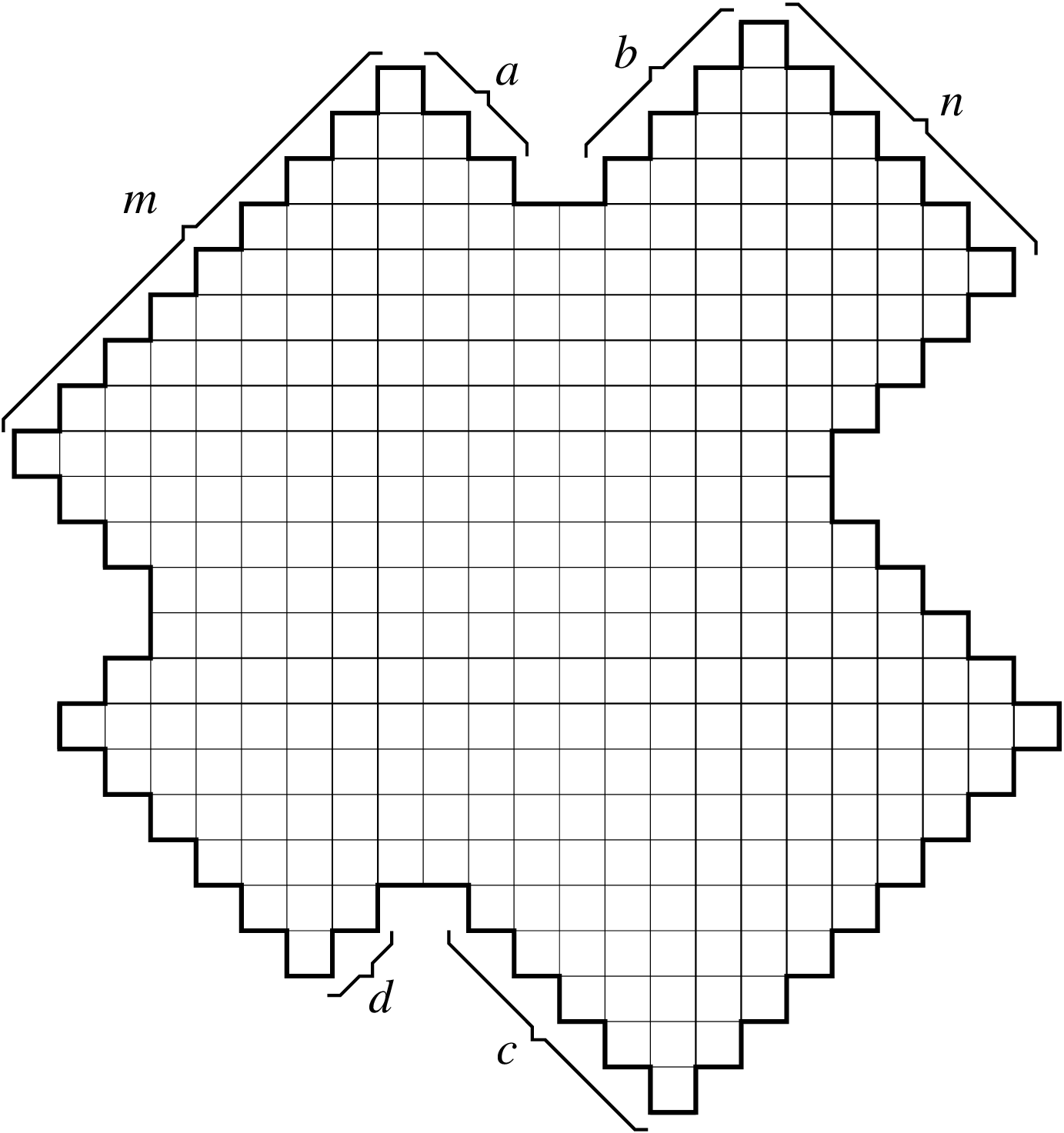}}
\hfill
{\includegraphics[width=0.43\textwidth]{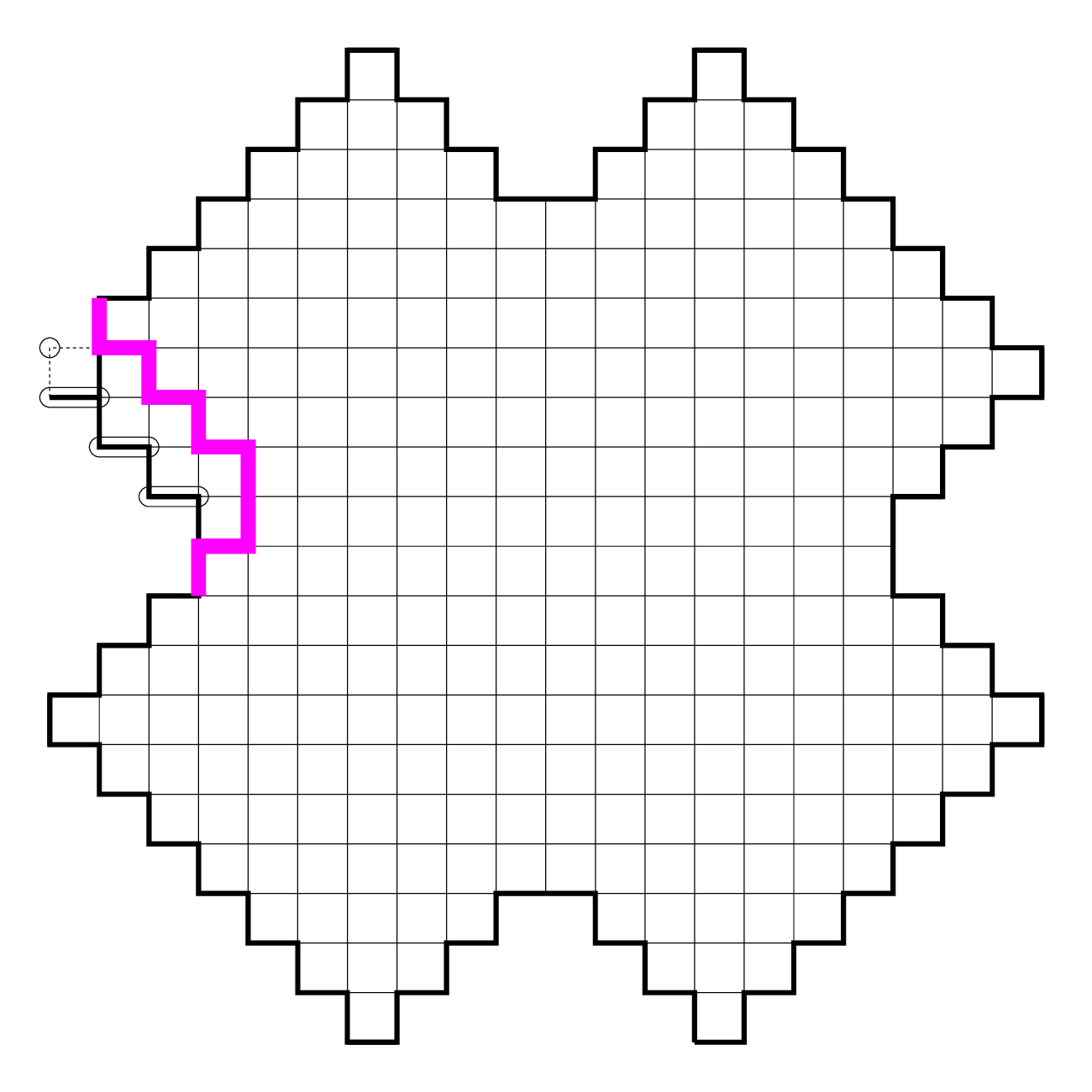}}
}
\vskip-0.1in
\caption{{\it Left.} The cruciform graph $C_{m,n}^{a,b,c,d}$ for $m=9$, $n=6$, $a=3$, $b=4$, $c=5$ and $d=2$. {\it Right.} The graph $\dot{C}_{2n+1,2n+1}^{n,n,n,n}$ has the same number of perfect matchings as the near-cruciform graph $D_{2n+1,2n+1}^{n,n,n,n}$. Here and throughout the paper, circled edges indicate forced edges (i.e., edges that are contained in every perfect matching).}
\label{fbaa}
\end{figure}

The {\it cruciform graph} $C_{m,n}^{a,b,c,d}$ is the subgraph of the square grid described on the left in Figure~\ref{fbaa} --- the four ``piers'' have lengths\footnote{ For the case when some of the pier lengths are negative we refer the reader to \cite[Figure 4]{Ciu2022}.} $a$, $b$, $c$, $d$ (clockwise from top left), the northwestern and southeastern piers (which have to be aligned) have widths $m$, and the other two (which are also aligned) have widths $n$. As shown in \cite{Ciu2022}, $C_{m,n}^{a,b,c,d}$ is balanced (i.e., has the same number of vertices in its two bipartition classes) if and only if $a+b+c+d=m+n-1$.

\begin{theo}[Ciucu,\cite{Ciu2022}]
\label{tba}
Let $C_{m,n}^{a,b,c,d}$ be a balanced cruciform graph. Then\footnote{ $\M(G)$ denotes the number of perfect matchings of the graph $G$.}
\begin{align}
\nonumber
\\[0pt]
\M(C_{m,n}^{a,b,c,d})&=2^{\left\{\frac14m(3m+1)+\frac14n(3n+1)-\frac12(a+c)(b+d)-\frac14(m-n)(a-b+c-d)\right\}}
\nonumber
\\[5pt]
&\ \ \ \ 
\times
\frac{\h(m+n+1)^2\h(m-a)\h(n-b)\h(m-c)\h(n-d)}{\h(n+a+1)\h(m+b+1)\h(n+c+1)\h(m+d+1)},
\label{ebaa}
\end{align}
where\footnote{ We set $\h(0)=1$ and $\h(n)=0$ for $n<0$.}
\begin{equation}
\h(n):=0!\,1!\cdots (n-1)!
\label{ebba}
\end{equation}

\end{theo}

The second result is a counterpart of the one above. It concerns a new family of graphs whose number of matchings is given by a simple product formula. We were led to it during the course of working out the simple proof of formula \eqref{eaa} which forms the subject of this paper.

Let $\dot{C}_{m,n}^{a,b,c,d}$ be the graph obtained from the cruciform graph $C_{m,n}^{a,b,c,d}$ by deleting the higher of the two leftmost vertices in the northwestern pier (see the picture on the left in Figure \ref{fbb}).
Let $D_{m,n}^{a,b,c,d}$ be the graph obtained from $\dot{C}_{m,n}^{a,b,c,d}$ after discarding all the forced edges in it (see the picture on the right in Figure \ref{fbaa}); we call it a {\it nearly-cruciform graph}. One readily checks that $D_{m,n}^{a,b,c,d}$ is balanced precisely if $a+b+c+d=m+n-2$.
By construction, $D_{m,n}^{a,b,c,d}$ and $\dot{C}_{m,n}^{a,b,c,d}$ have the same number of matchings.


For general values of the parameters, the  number of matchings of $D_{m,n}^{a,b,c,d}$ does not seem to be given by a simple product formula (indeed, in concrete examples, relatively large primes come up in its factorization). However, for $a=c$, the formula below holds. It is a strike of good luck that this is precisely the case we need for our proof of formula \eqref{eaa}.

\begin{theo}
\label{tbb}
    Let $D_{m,n}^{a,b,a,d}$ be a balanced nearly-cruciform graph. Then
    \begin{equation}\label{ebca}
    \begin{aligned}
        \M(D_{m,n}^{a,b,a,d})=&2^{\{n(n-2a-2)-3n(n-1)/2+(m-a)(m-a-1)+(n+a+1)(n+a)\}}\frac{(m+b-1)!!}{(n-d-1)!!}\\
        &\times\frac{\h(m+n+1)\h(m+n)\h(m-a)\h(m-a-1)\h(n-b)\h(n-d)}{\h(n+a+1)^2\h(m+b+1)\h(m+d+1)}.
    \end{aligned}
    \end{equation}    
\end{theo}

\medskip
The proof of Theorem \ref{tbb} is presented in Section 3.

\begin{figure}[t]
\centerline{
{\includegraphics[width=0.47\textwidth]{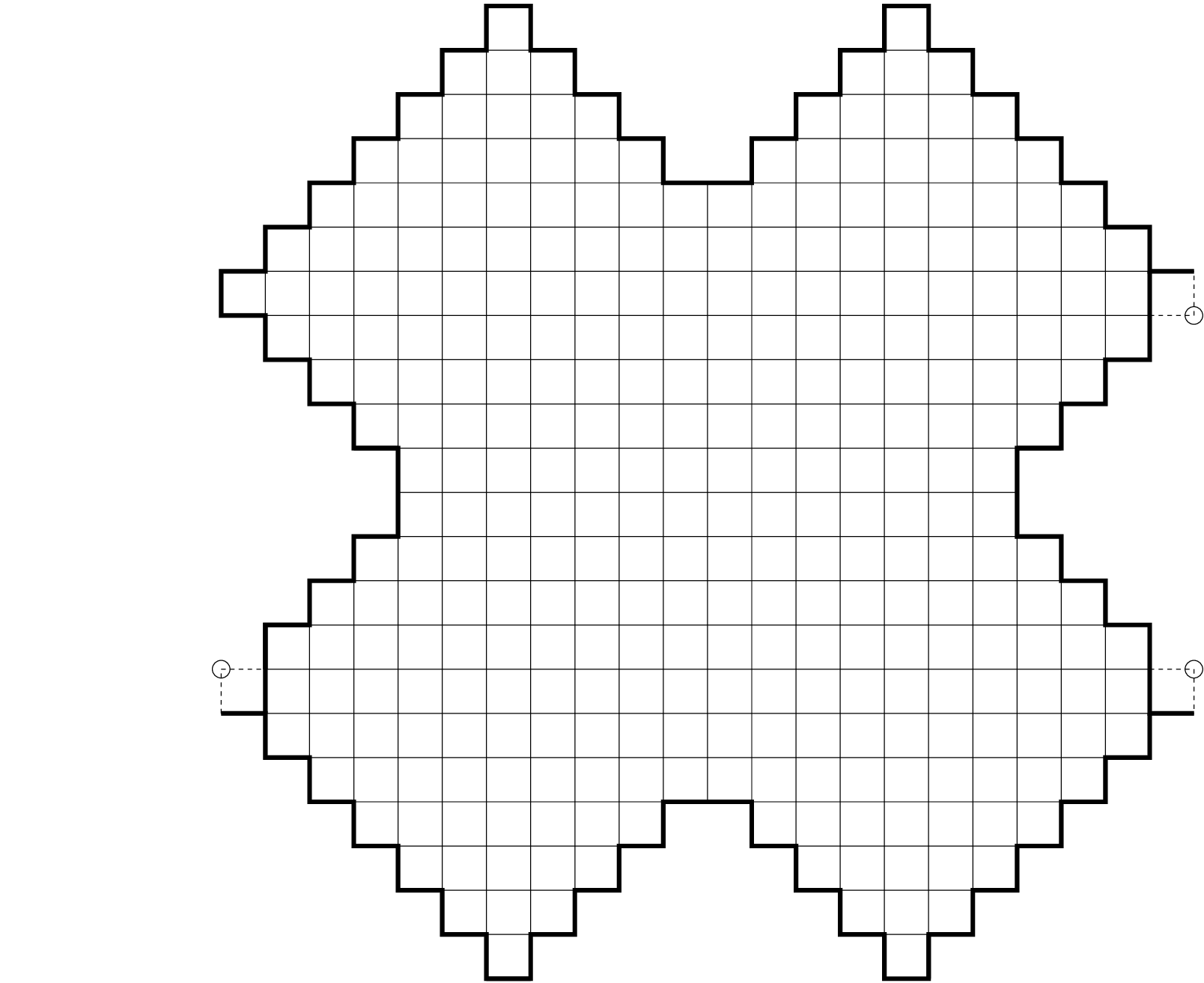}}
\hfill
{\includegraphics[width=0.40\textwidth]{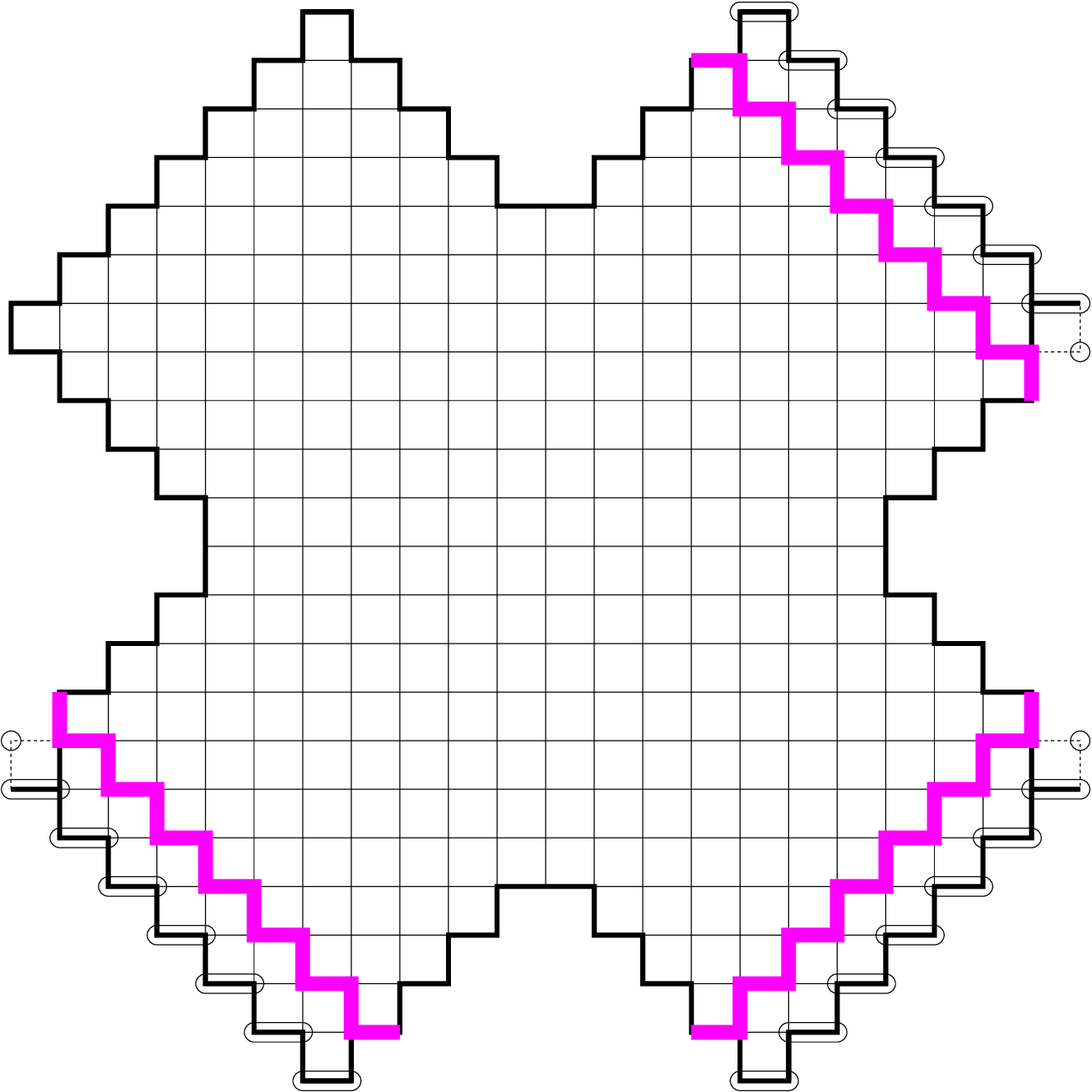}}
}
\vskip0.5in
\centerline{
{\includegraphics[width=0.47\textwidth]{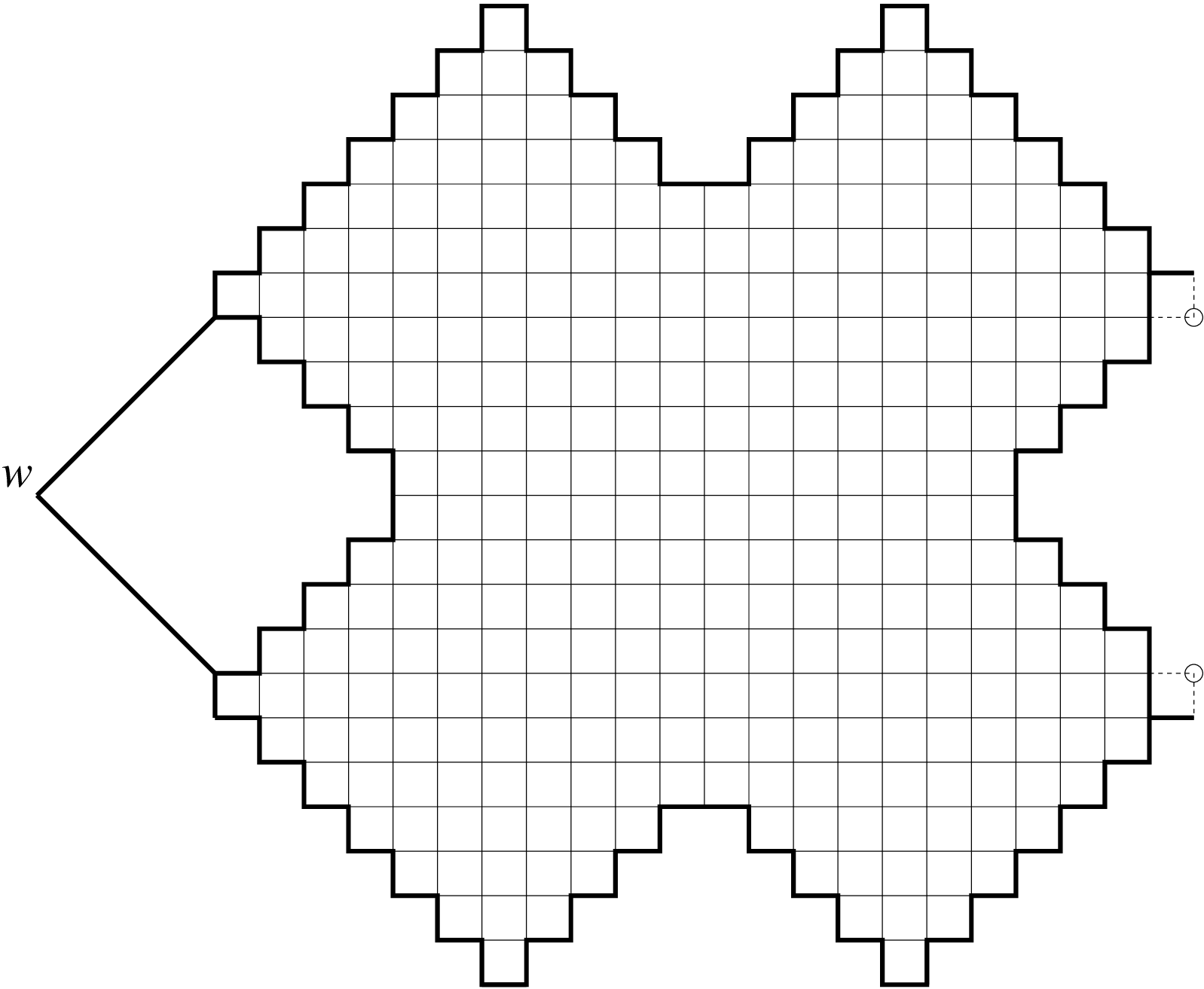}}
\hfill
{\includegraphics[width=0.47\textwidth]{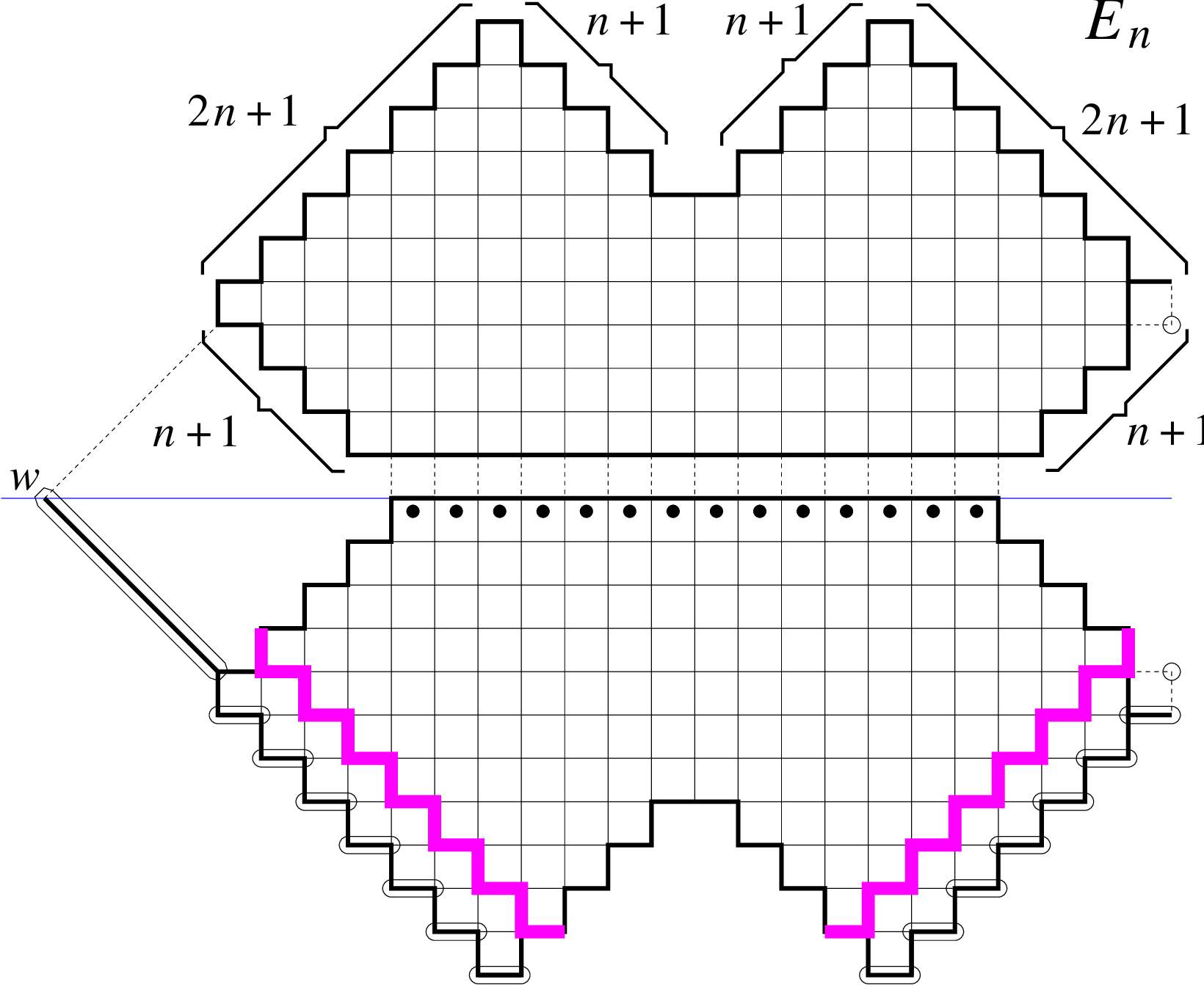}}
}
\caption{ Obtaining equation \eqref{eba} (here $n=3$). {\it Top left.} $\dddot{C}_{2n+1,2n+1}^{n+1,n+1,n+1,n+1}$ is obtained from the cruciform graph $C_{2n+1,2n+1}^{n+1,n+1,n+1,n+1}$ by removing the three indicated boundary vertices.
{\it Top right.} Due to forced edges, $\dddot{C}_{2n+1,2n+1}^{n+1,n+1,n+1,n+1}$ has the same number of matchings as the cruciform graph $C_{2n+1,2n+1}^{n+1,n,n,n}$ (the portion of the grid graph contained within the thick contour that has three magenta portions).
\newline
{\it Bottom left.} A graph with twice as many matchings as $\dddot{C}_{2n+1,2n+1}^{n+1,n+1,n+1,n+1}$.
{\it Bottom right.} Applying the factorization theorem to the graph with twice as many matchings as $\dddot{C}_{2n+1,2n+1}^{n+1,n+1,n+1,n+1}$, so as to produce the graph $E_n$ on top and (after the removal of the forced edges) the graph $F_n$ on the bottom.
}
\label{fba}
\end{figure}

\begin{figure}[t]
\centerline{
{\includegraphics[width=0.43\textwidth]{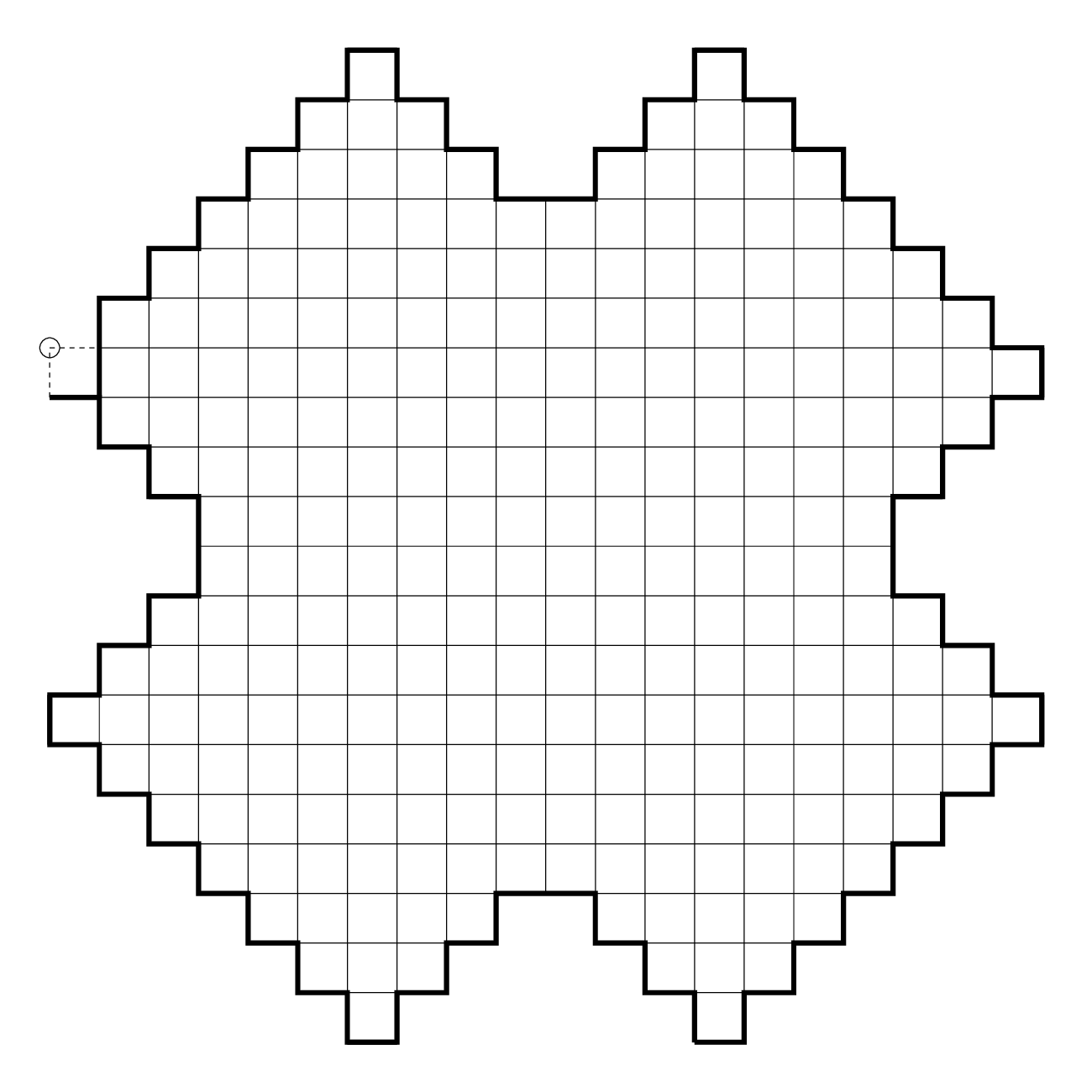}}
\hfill
{\includegraphics[width=0.52\textwidth]{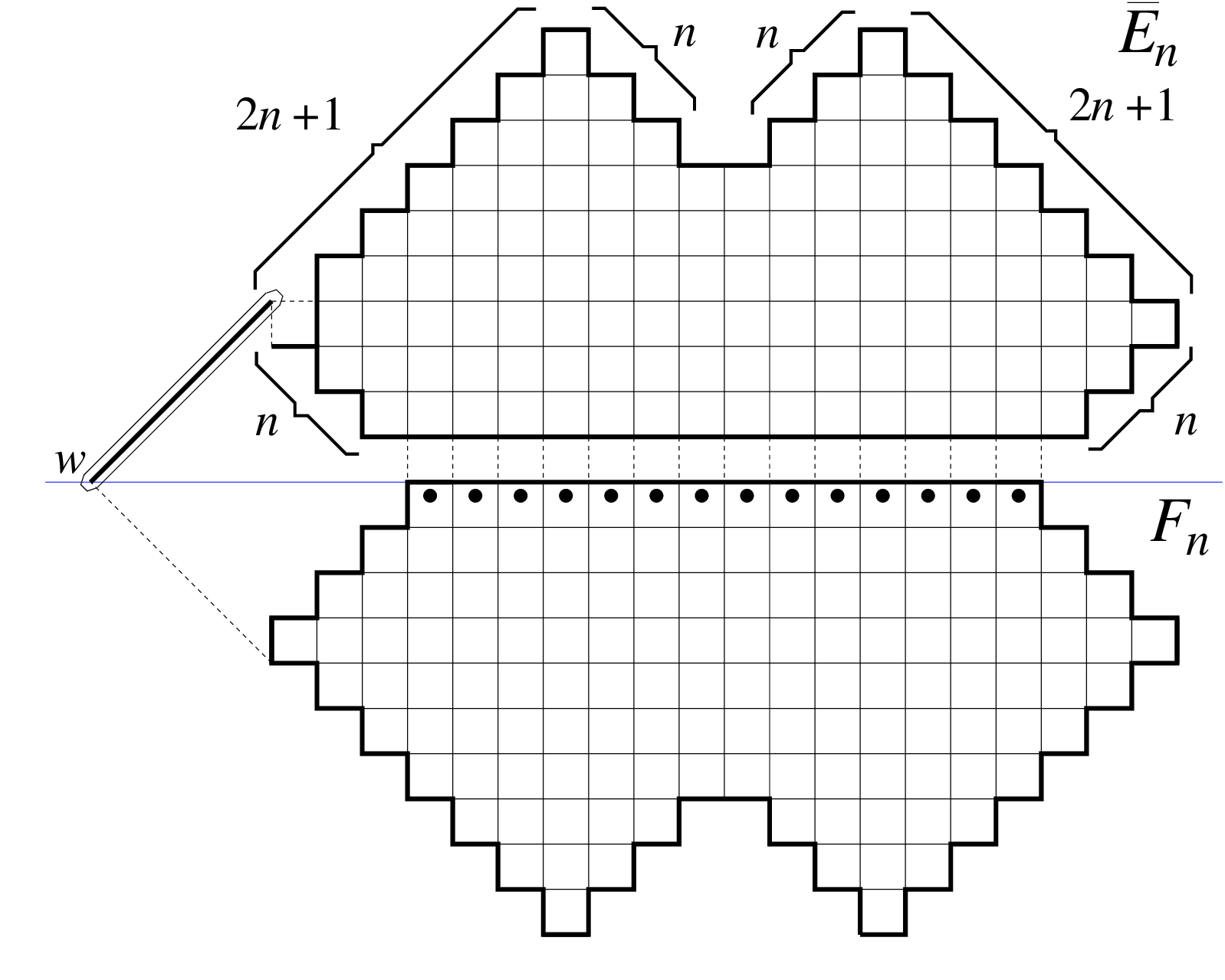}}
}
\vskip-0.1in
\caption{{\it Left.} The graph $\dot{C}_{2n+1,2n+1}^{n,n,n,n}$  for $n=3$.
{\it Right.} Applying the factorization theorem to a graph with twice as many perfect matchings as $\dot{C}_{2n+1,2n+1}^{n,n,n,n}$ produces the graph $\overline{E}_n$ on top and the same graph $F_n$ on the bottom as on the bottom right in Figure~\ref{fba}.}
\vskip-0.1in
\label{fbb}
\end{figure}

\medskip
{\it Proof of formula \eqref{eaa}.}
Let $\dddot{C}_{2n+1,2n+1}^{n+1,n+1,n+1,n+1}$ be the graph obtained from the cruciform graph ${C}_{2n+1,2n+1}^{n+1,n+1,n+1,n+1}$ by removing three boundary vertices as indicated on the top left of Figure \ref{fba}.
This graph is bipartite, and almost symmetric with respect to a horizontal symmetry axis (the only exception is on the left boundary). So we are close to being able to apply the factorization theorem of \cite{Ciu1997} to it, and in fact we can do it, using the following trick (this will be used several times throughout the paper). 

Consider the graph $G$ obtained from $\dddot{C}_{2n+1,2n+1}^{n+1,n+1,n+1,n+1}$ by including back in it the deleted vertex from the southwestern pier and its two incident edges, as well as a new vertex $w$ and two new edges incident to $w$, as shown on the bottom left in Figure \ref{fba}. Due to the horizontal symmetry of $G$ and to the fact that one of the two new edges is contained in each perfect matching of $G$, it follows that $G$  has exactly twice as many matchings as $\dddot{C}_{2n+1,2n+1}^{n+1,n+1,n+1,n+1}$.

Furthermore,~$G$ is a symmetric planar bipartite graph, and thus we can apply to it the factorization theorem \cite[Theorem 2.1]{Ciu1997}. This theorem states that $M(G)=2^k\M(G^+)\M(G^-)$, where $G^+$ and $G^-$ are ``halves'' obtained from $G$ by performing the following operations:

\medskip
$(i)$ label the vertices on the horizontal\footnote{ We can assume without loss of generality that the symmetry axis is horizontal (otherwise we rotate the graph).} symmetry axis $\ell$ from left to right $a_1,b_1,\dotsc,a_k,b_k$

$(ii)$ color the vertices on $\ell$ white and black according to the two color classes of $G$

$(iii)$ remove from $G$ all the edges incident from above to each white $a_i$ and each black $b_i$

$(iv)$ remove from $G$ all the edges incident from below to each black $a_i$ and each white $b_i$

$(v)$ place weight\footnote{ For a graph $G$ with weights on its edges, $\M(G)$ is the total weight of its perfect matchings. The weight of a perfect matching is the product of the weights of its constituent edges (unless another weight is specified, each edge of a graph carries weight 1).} $1/2$ on each edge of $G$ along $\ell$

\medskip
Therefore, using also the fact that ${C}_{2n+1,2n+1}^{n+1,n,n,n}$ and $\dddot{C}_{2n+1,2n+1}^{n+1,n+1,n+1,n+1}$ have the same number of matchings (this is due to forced edges; see the top right picture in Figure \ref{fba}), we obtain that
\begin{equation}
2\M({C}_{2n+1,2n+1}^{n+1,n,n,n})=2\M(\dddot{C}_{2n+1,2n+1}^{n+1,n+1,n+1,n+1})=\M(G)=2^{2n+2}\M(E_n)\M(F_n),
\label{eba}
\end{equation}
where $E_n$ is the graph $G^+$ above $\ell$ in the bottom right picture of Figure \ref{fba}, and $F_n$ is the graph obtained from $G^-$ after removing the forced edges.

Next, we apply the same reasoning to the graph $\dot{C}_{2n+1,2n+1}^{n,n,n,n}$ (pictured on the left in Figure~\ref{fbb}). More precisely, we augment it to a graph $G_1$ by restoring the deleted vertex and including a new vertex and two new edges as shown in the picture on the right in Figure \ref{fbb}, and we apply the factorization theorem to $G_1$. We obtain
\begin{equation}
2\M({D}_{2n+1,2n+1}^{n,n,n,n})=2\M(\dot{C}_{2n+1,2n+1}^{n,n,n,n})=\M(G_1)=2^{2n+2}\M(\overline{E}_n)\M(F_n),
\label{ebb}
\end{equation}
where $\overline{E}_n$ is the graph above $\ell$ in the picture on the right in Figure \ref{fba}, and $F_n$, the graph below $\ell$ resulting from the factorization theorem is --- very conveniently! --- precisely the same as the one in equation \eqref{eba} (compare the right picture in Figure \ref{fbb} and the bottom right picture in Figure \ref{fba}). This invites us to divide side by side equations \eqref{eba} and \eqref{ebb}, to obtain
\begin{equation}
\frac{\M(C^{n+1,n,n,n}_{2n+1,2n+1})}{\M(D^{n,n,n,n}_{2n+1,2n+1})}=\frac{\M(E_n)}{\M(\overline{E}_{n})}.
\label{ebc}
\end{equation}
To anticipate how our proof will come to completion, we can say that equation \eqref{ebc} and its proof represents one of three steps which, when linked together, will imply that
\begin{equation}
    \frac{\M(\mathcal{T}_{n+1})}{\M(\mathcal{T}_n)}=\frac{\M(C^{n+1,n,n,n}_{2n+1,2n+1})}{2\M(D^{n,n,n,n}_{2n+1,2n+1})}.
\label{ebd}
\end{equation}
Since the quantities on the right hand side above are given by the formulas from Theorems \ref{tba} and \ref{tbb}, formula \eqref{eaa} will follow by simply verifying that the resulting product on the right hand side above agrees with the expression for the left hand side above resulting from \eqref{eaa}.

\begin{figure}[t]
\centerline{
{\includegraphics[width=0.47\textwidth]{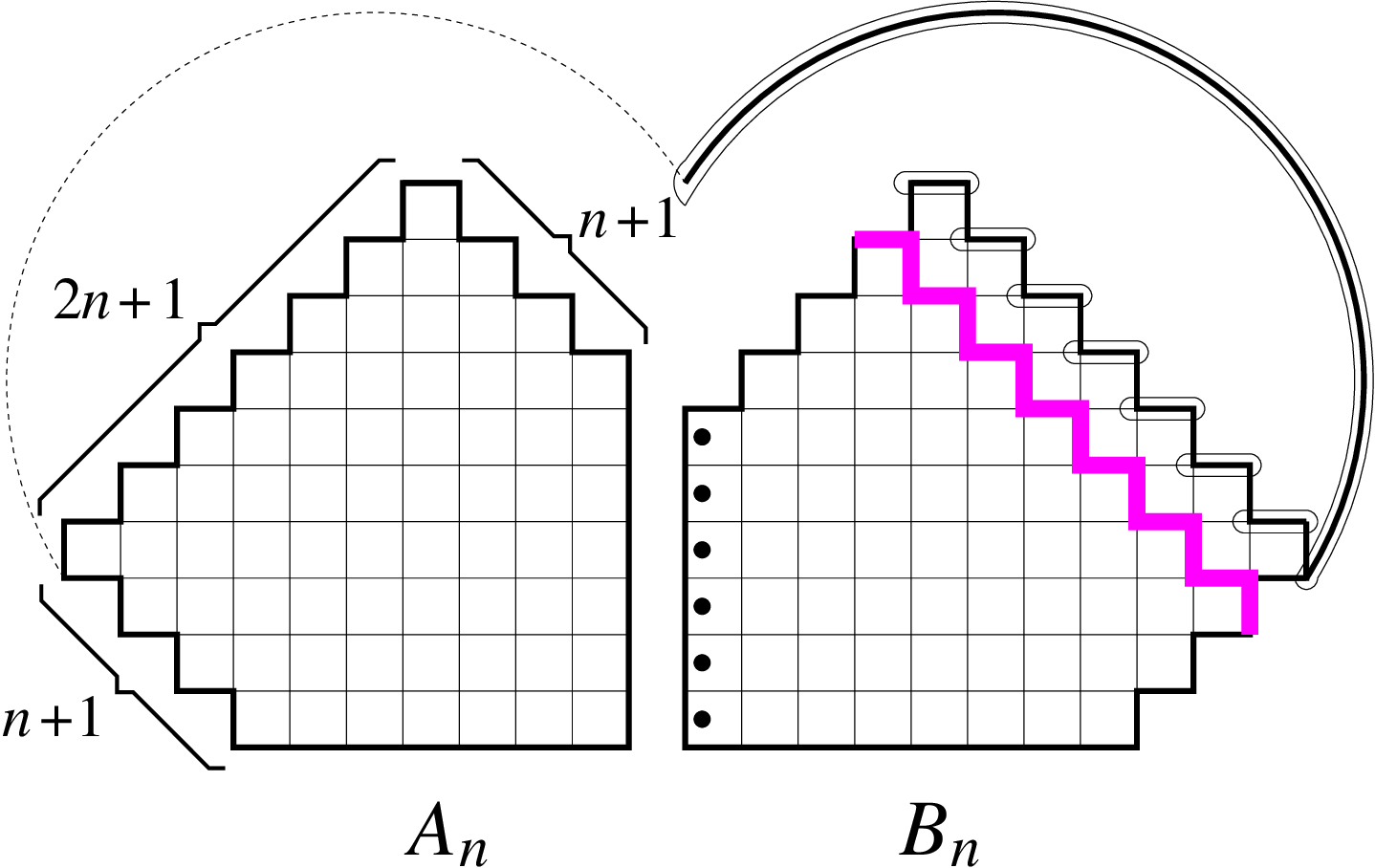}}
\hfill
{\includegraphics[width=0.425\textwidth]{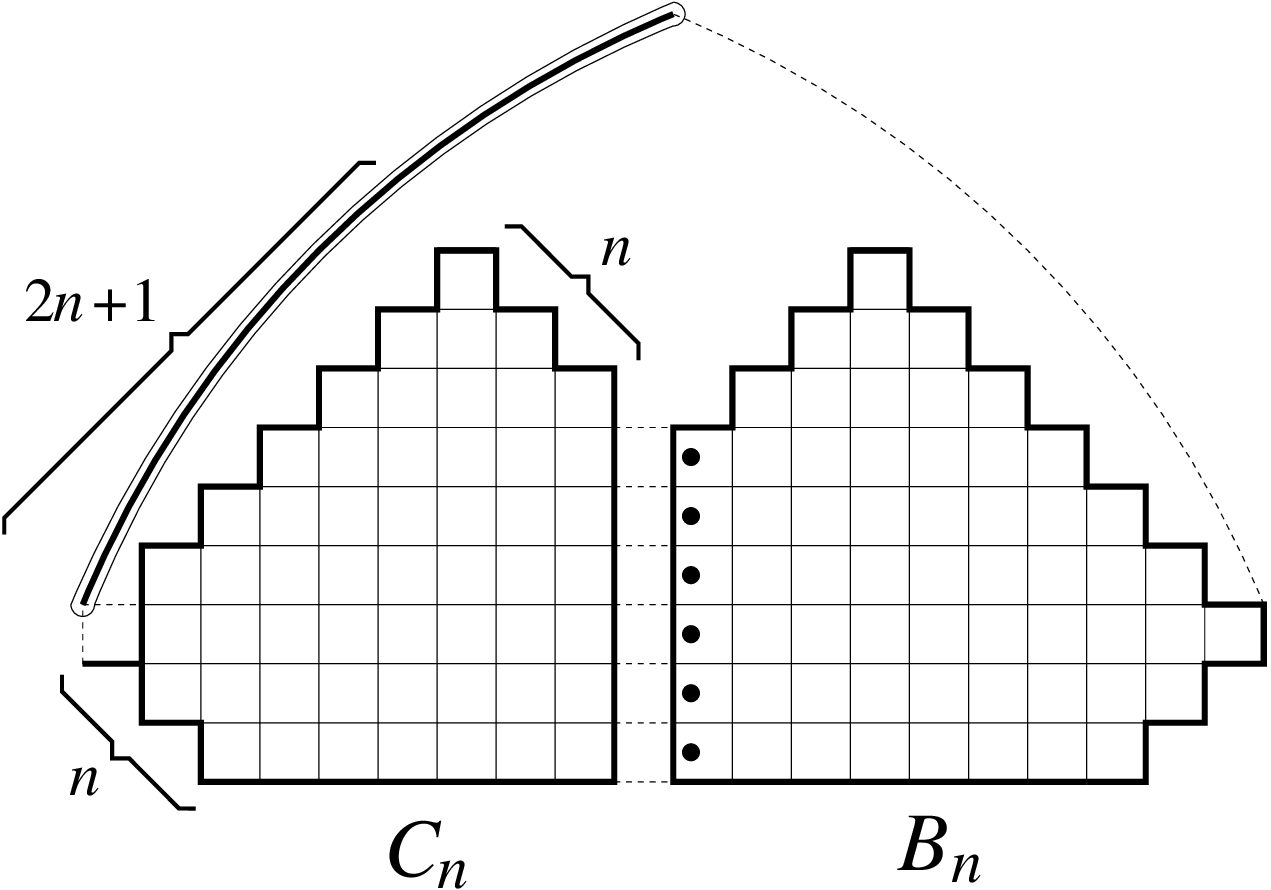}}
}
\vskip-0.1in
\caption{{\it Left.} Obtaining equation \eqref{ebe} (here $n=3$). When applying the factorization theorem to a graph having twice as many matchings as $E_n$, the resulting $G^+$ is the graph $A_n$, while the resulting $G^-$ is the graph $B_n$.
{\it Right.} Obtaining equation \eqref{ebf} (here $n=3$). When applying the factorization theorem to a graph having twice as many matchings as $\overline{E}_n$, the resulting $G^+$ is the graph $C_n$, while the resulting $G^-$ is the same graph $B_n$ as in the picture on the left.}
\vskip-0.1in
\label{fbc}
\end{figure}

To continue with our proof, in the second step we apply the factorization theorem to the graphs $E_n$ and $\overline{E}_n$. For $E_n$, this is shown on the left in Figure \ref{fbc}. Following the same approach as before, we consider a symmetric graph $G_2$ that has twice as many matchings as $E_n$, and apply the factorization theorem to it. One readily verifies that applying the steps $(i)$--$(v)$, the resulting graphs $G_2^+$ and $G_2^-$ are the graphs $A_n$ and $B_n$ shown on the left in Figure \ref{fbc}. There are $2n+2$ vertices on the symmetry axis, so the exponent of 2 is $n+1$. We obtain
\begin{equation}
2\M(E_n)=\M(G_2)=2^{n+1}\M(A_n)\M(B_n).
\label{ebe}
\end{equation}
The same procedure, when applied to the graph $\overline{E}_n$, is illustrated in the picture on the right in Figure \ref{fbc}. The $G^+$ graph resulting from applying the factorization theorem is now the graph $C_n$ described in that picture, while the $G^-$ graph is exactly the same graph $B_n$ that appears on the left in Figure \ref{fbc} and in equation \eqref{ebe}. This gives
\begin{equation}
2\M(\overline{E}_n)=2^{n+1}\M(C_n)\M(B_n).
\label{ebf}
\end{equation}
Combining equations \eqref{ebe} and \eqref{ebf} we get
\begin{equation}
\frac{\M({E}_n)}{\M(\overline{E}_n)}=\frac{\M({A}_n)}{\M({C}_n)}.
\label{ebg}
\end{equation}
\begin{figure}[t]
\centerline{
\hfill
{\includegraphics[width=0.25\textwidth]{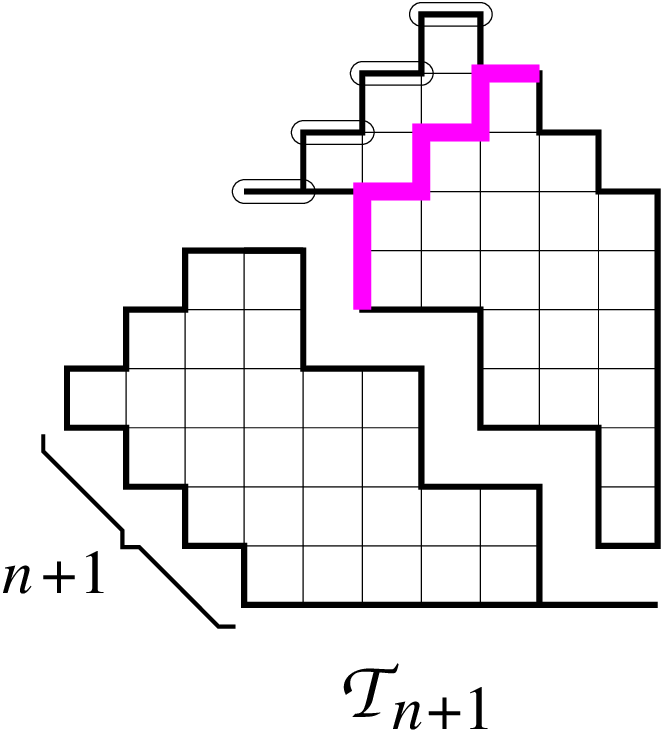}}
\hfill
{\includegraphics[width=0.22\textwidth]{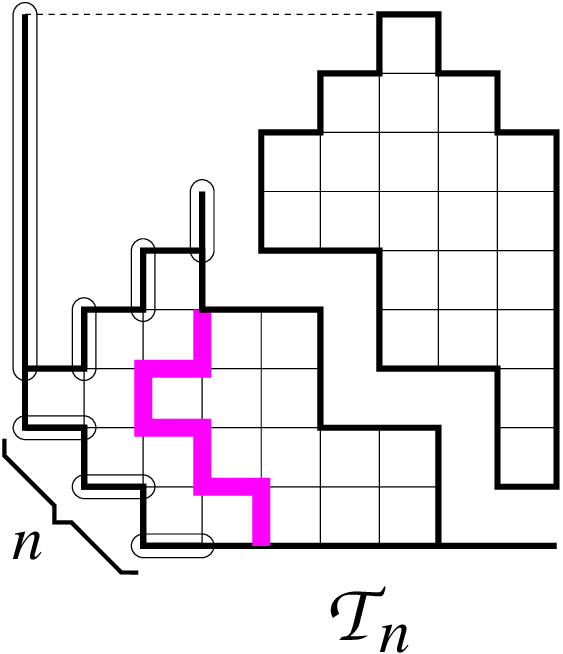}}
\hfill
}
\vskip-0.1in
\caption{{\it Left.} Obtaining equation \eqref{ebh} (here $n=3$). When applying the factorization theorem to the graph $A_n$, the resulting $G^+$ is the Aztec triangle graph ${\mathcal T}_{n+1}$, while the resulting $G^-$ is the graph $S_n$.
{\it Right.} Obtaining equation \eqref{ebi} (here $n=3$). When applying the factorization theorem to a graph having twice as many matchings as $C_n$, the resulting $G^+$ is the Aztec triangle graph ${\mathcal T}_{n}$, while the resulting $G^-$ is precisely the same graph $S_n$ as in the picture on the left.}
\vskip-0.1in
\label{fbd}
\end{figure}

We are now one step away from obtaining \eqref{ebd}. Apply the factorization theorem to the graph $A_n$ (in this case we can do it directly; see the picture on the left in Figure \ref{fbd}). This time the symmetry axis is along a diagonal of the square grid, so all vertices on it have the same color. Due to this, the procedure $(i)$-$(v)$ prescribes to alternate between deleting the edges incident from below and deleting those incident from above, as we go successively through the vertices on $\ell$. As it is clear from Figure~\ref{fbd}, the resulting $G^+$ graph is precisely the Aztec triangle graph ${\mathcal T}_{n+1}$ (compare with the picture on the right in Figure \ref{faa}); denote the resulting $G^-$ graph by $S_n$. The number of vertices on the symmetry axis is $2n+2$. Thus we obtain
\begin{equation}
\M(A_n)=2^{n+1}\M({\mathcal T}_{n+1})\M(S_n).
\label{ebh}
\end{equation}

Finally, apply our procedure to the graph $C_n$. As usual, consider a symmetric graph that has twice as many matchings as $C_n$, and apply the factorization theorem to it. The resulting $G^+$ graph is the Aztec triangle graph ${\mathcal T}_{n}$, while the $G^-$ graph is precisely the same graph $S_n$ as on the left in Figure \ref{fbd} and in equation \eqref{ebh}. Since the number of vertices on the symmetry axis is still $2n+2$, this gives
\begin{equation}
2\M(C_n)=2^{n+1}\M({\mathcal T}_n)\M(S_n).
\label{ebi}
\end{equation}
Combining equations \eqref{ebh} and \eqref{ebi} we obtain
\begin{equation}
\frac{\M({A}_n)}{\M({C}_n)}=2\frac{\M({\mathcal T}_{n+1})}{\M({\mathcal T}_n)},
\label{ebj}
\end{equation}
which, together with \eqref{ebc} and \eqref{ebg} prove equation \eqref{ebd}.

Using Theorems \ref{tba} and \ref{tbb} and some straightforward manipulations we obtain
\begin{equation}
\begin{aligned}
    \frac{\M(\mathcal{T}_{n+1})}{\M(\mathcal{T}_n)}=\frac{\M(C^{n+1,n,n,n}_{2n+1,2n+1})}{2\M(\dot{C}^{n,n,n,n}_{2n+1,2n+1})}&=2^n\frac{n!!}{(3n)!!}\frac{(4n+2)!}{(3n+2)!}\\
    &=\Bigg[2^{(n+1)n/2}\prod_{i=0}^{n}\frac{(4i+2)!}{(n+2i+2)!}\Bigg]/\Bigg[2^{n(n-1)/2}\prod_{i=0}^{n-1}\frac{(4i+2)!}{(n+2i+1)!}\Bigg].
\end{aligned}
\label{ebk}
\end{equation}
In other words, the ratio between $\M(\mathcal{T}_{n+1})$ and $\M(\mathcal{T}_n)$ is equal to the ratio of the expressions for these that follow from the claimed formula \eqref{eaa}. 
Since the two sides of \eqref{eaa} clearly agree for $n=1$ (both sides being equal to 1), it follows that they agree for all $n$. 
$\hfill\square$

\section{Proof of Theorem \ref{tbb}}

The {\it Aztec rectangle} graph $AR_{m,n}$ is the graph whose vertices are the white squares of a $(2m+1)\times(2n+1)$ chessboard with black corners, two vertices being connected by an edge precisely if the corresponding white squares share a corner. The {\it trimmed Aztec rectangle} $\overline{AR}_{m,n}$ is the graph obtained from $AR_{m,n}$ by removing its $n$ bottom-most vertices. $\overline{AR}_{m,n}$ is bipartite, and has $n+1$ vertices on the bottom; one readily sees that for $m\leq n$, $(1)$ it has $m$ more vertices of the majority color, and $(2)$ the $n+1$ bottom vertices have the majority color.

Let the set $T=\{t_1,\cdots,t_m\}\subseteq[n+1]$ have its elements sorted in increasing order, and consider the graph obtained from the trimmed Aztec rectangle $\overline{AR}_{m,n}$ by deleting its $t_1$-st, $t_2$-nd, \dots, $t_m$-th (from left to right) vertex on the bottom; denote it by $\overline{AR}_{m,n}(T)$ (two instances are illustrated in the picture on the right in Figure \ref{fcc}).

We will employ the following result, due independently to Elkies, Kuperberg, Larsen and Propp \cite{EKLP} and to Helfgott and Gessel \cite[Lemma 3]{GH}. 

\begin{lem}
Given positive integers $m$ and $n$ with $m\leq n$ and a set $T=\{t_1,\dotsc,t_m\}\subseteq[n+1]$ with elements sorted in increasing order, the number of matchings of the graph $\overline{AR}_{m,n}(T)$ is given by
\begin{equation}
    \M(\overline{AR}_{m,n}(T))=\frac{2^{m(m-1)/2}}{\h(m)}\prod_{1\leq i<j\leq m}(t_j-t_i).
\label{eca}
\end{equation}    
\label{tca}
\end{lem}

{\it Proof of Theorem $\ref{tbb}$.} The proof follows by the same arguments that proved Theorem \ref{tba} in \cite{Ciu2022}: one applies successively the complementation theorem of \cite{Ciu1998} to the given nearly-cruciform graph, until it is gradually transformed into a graph whose number of matchings we can find easily. The only difference is that in \cite{Ciu2022} these successive transformations led to an Aztec rectangle graph with two intrusions (whose matchings were enumerated by Krattenthaler in \cite{Kra2000}), while in this case the final graph turns out to be a doubly-intruded Aztec rectangle with one vertex removed from a certain position on the boundary. Using the additional symmetry assumption in the statement of Theorem \ref{tbb} (namely, that the northwestern and the southeastern piers have the same length), we are able to deduce the number of matchings of the latter from Lemma \ref{tca}, using the factorization theorem.

The special case of the complementation theorem needed here is described in detail in the first two pages of \cite[Section 3]{Ciu2022}. In a nutshell, it states that the number of matchings of certain subgraphs $H$ of a cellular graph\footnote{A finite subgraph $G$ of the grid graph $\Z^2$ is called cellular if its set of edges can be partitioned into 4-cycles.} $G$ satisfies
\begin{equation}
\M(H)=2^t\M(H'),
\label{ecb}
\end{equation}
where the exponent $t$ has a simple explicit description, and the graph $H'$ (called the complement of $H$) is obtained from $H$ by a simple, explicit construction.

\begin{figure}[t]
\centerline{
\hfill
{\includegraphics[width=0.34\textwidth]{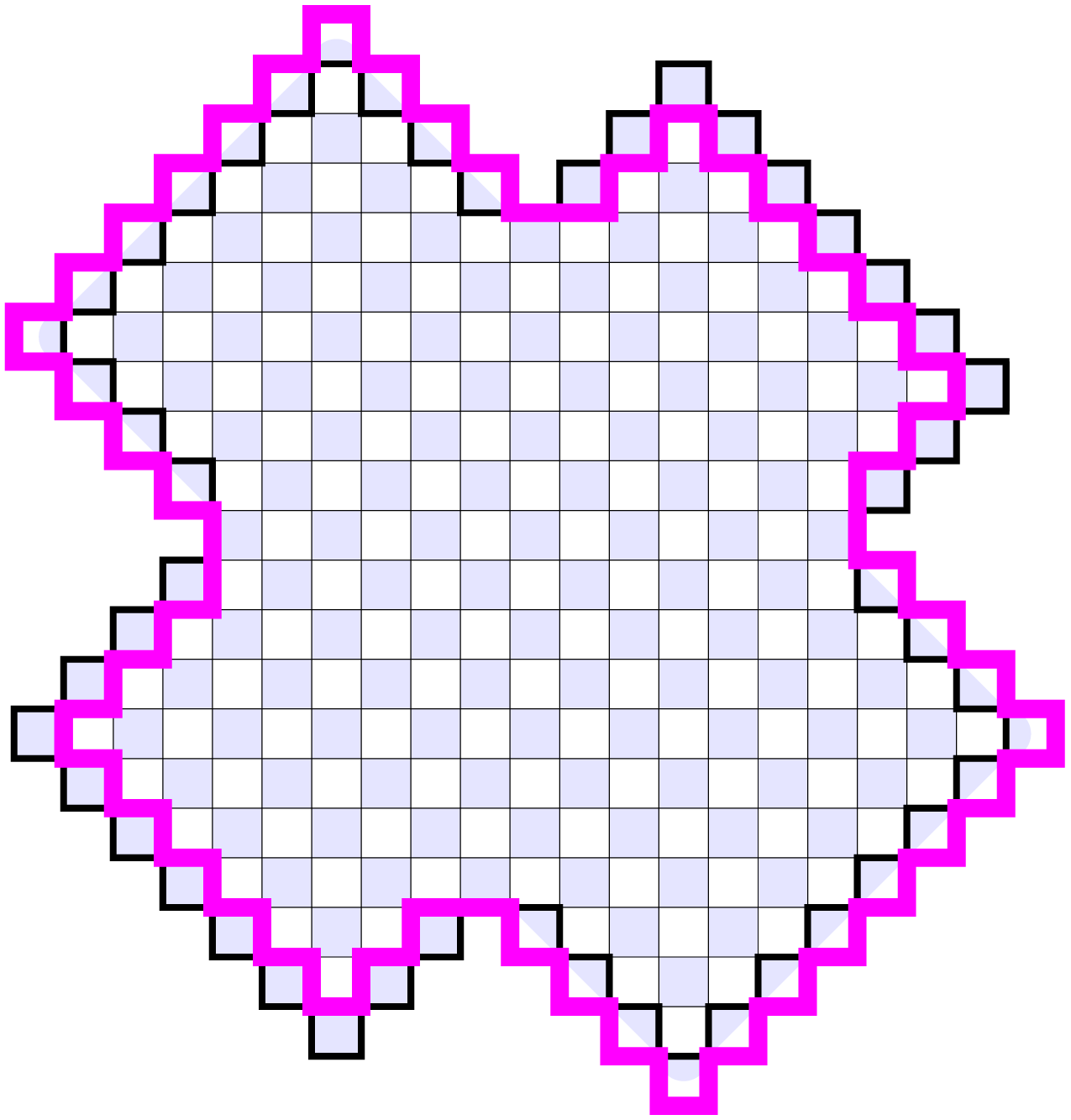}}
\hfill
{\includegraphics[width=0.38\textwidth]{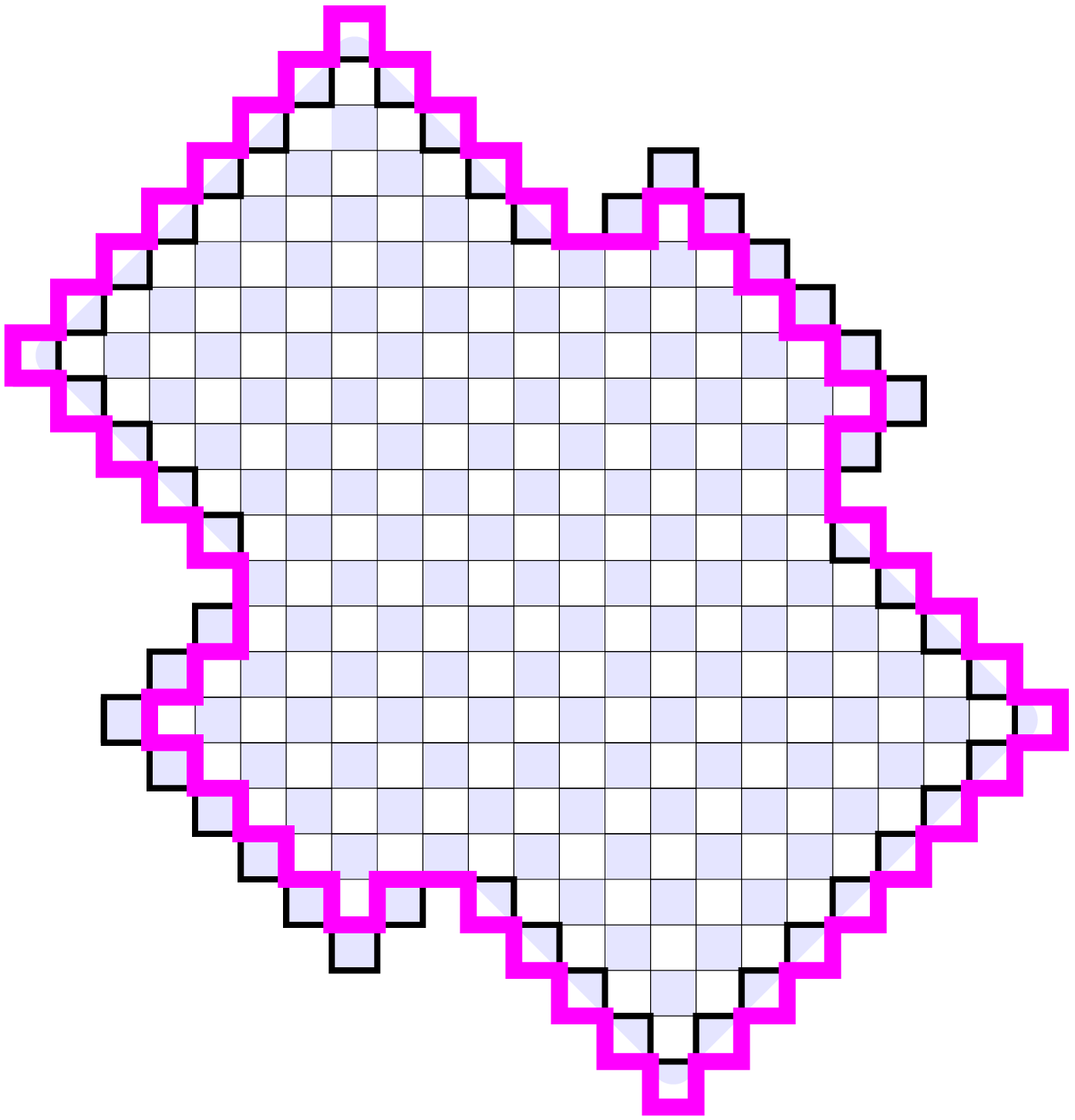}}
\hfill
}
\vskip0.1in
\centerline{
\hfill
{\includegraphics[width=0.40\textwidth]{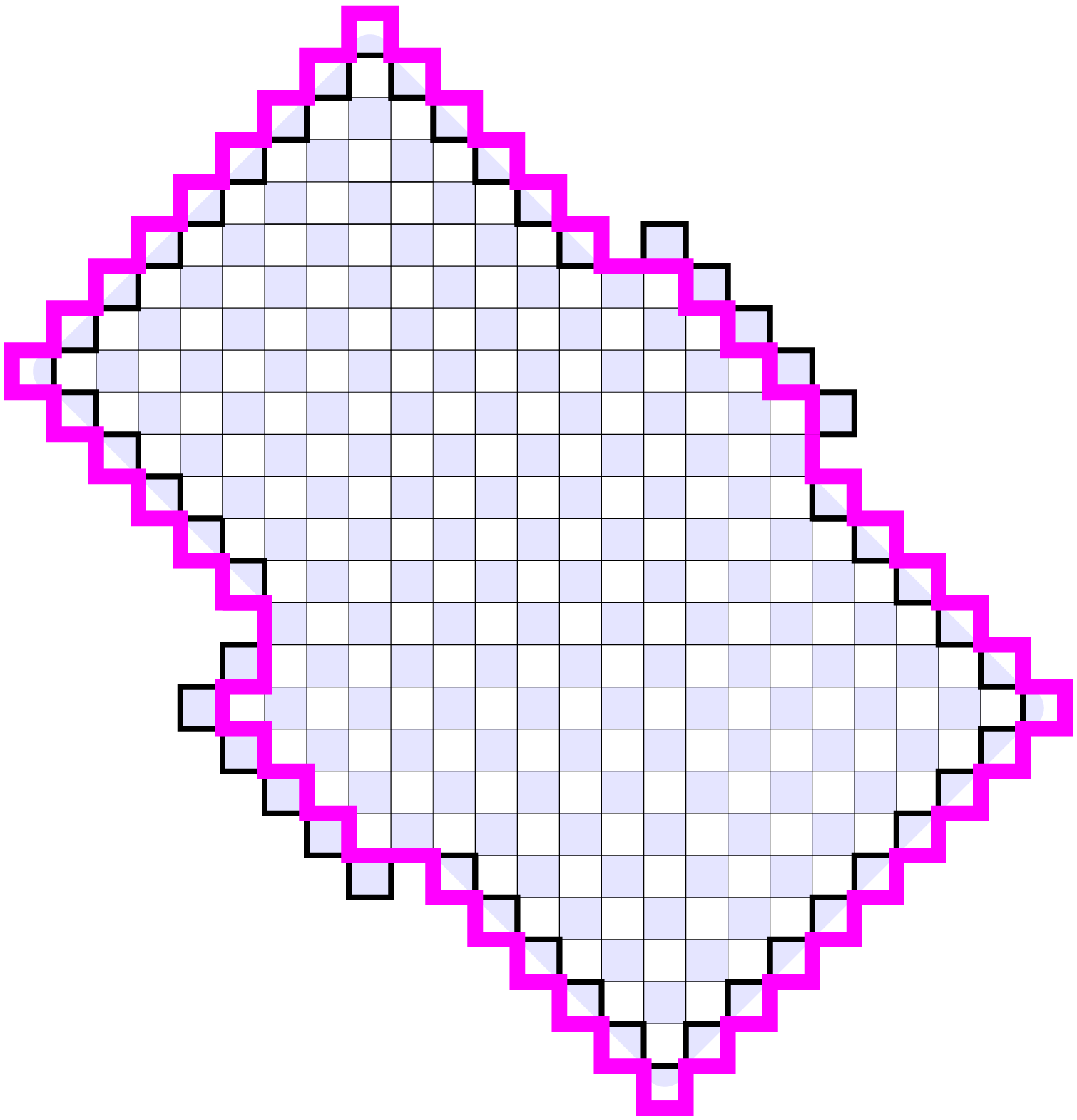}}
\hfill
{\includegraphics[width=0.43\textwidth]{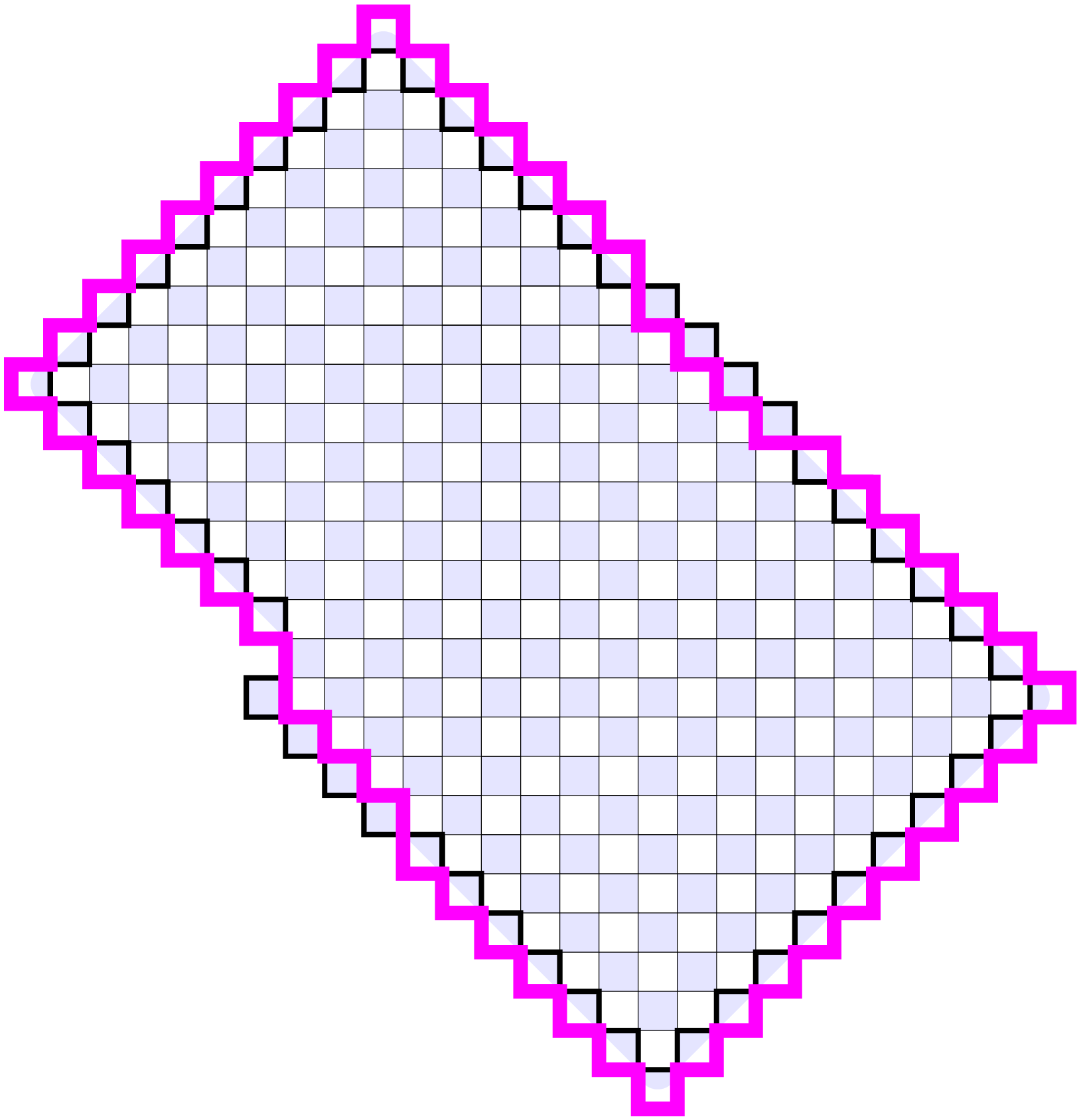}}
\hfill
}
\vskip0.1in
\caption{ Successive applications of the complementation theorem to the nearly-cruciform graphs: $D_{7,7}^{3,3,3,3}$ yields $D_{7+1,7-1}^{3+1,3-1,3+1,3-1}$, i.e.\ $D_{8,6}^{4,2,4,2}$ (top left); $D_{8,6}^{4,2,4,2}$  yields $D_{8+1,6-1}^{4+1,2-1,4+1,2-1}$, i.e.\ $D_{9,5}^{5,1,5,1}$ (top right); $D_{9,5}^{5,1,5,1}$  yields $D_{9+1,5-1}^{5+1,1-1,5+1,1-1}$, i.e.\ $D_{10,4}^{6,0,6,0}$ (bottom left); $D_{10,4}^{6,0,6,0}$  yields $D_{10+1,4-1}^{6+1,0-1,6+1,0-1}$, i.e.\ $D_{11,3}^{7,-1,7,-1}$ (bottom right).
}
\vskip-0.1in
\label{fca}
\end{figure}

\begin{figure}[t]
\centerline{
\hfill
{\includegraphics[width=0.42\textwidth]{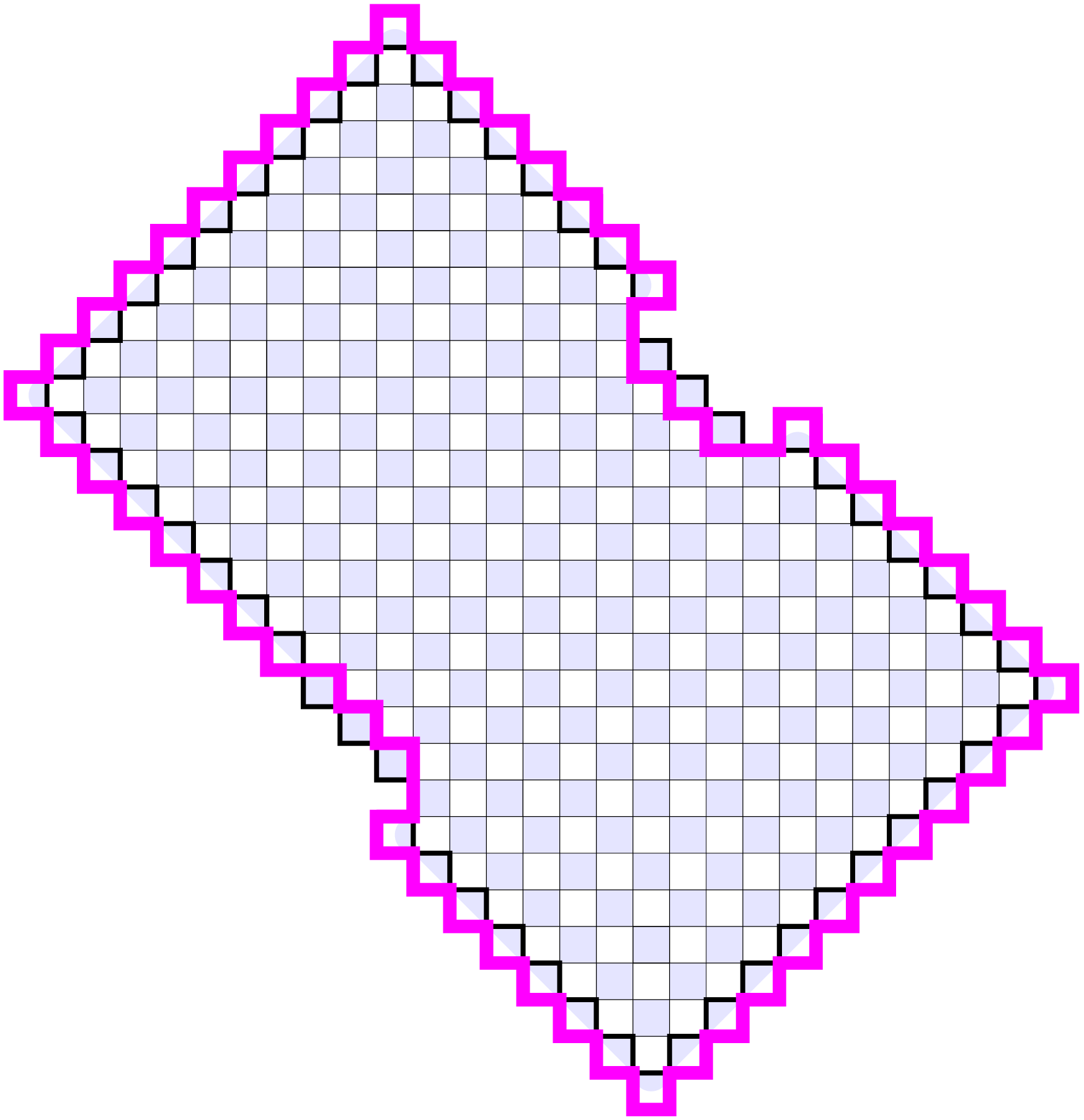}}
\hfill
{\includegraphics[width=0.445\textwidth]{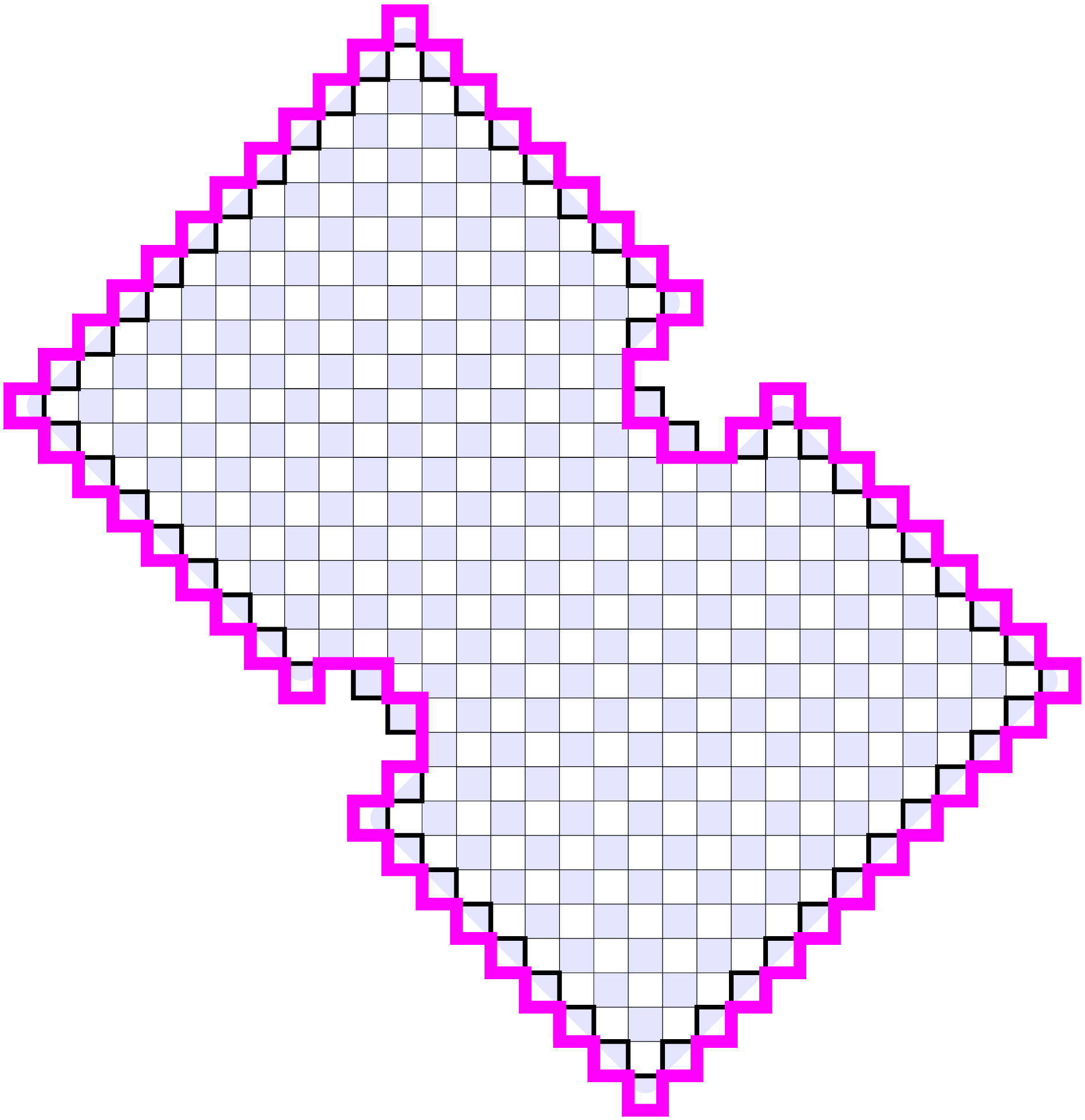}}
\hfill
}
\vskip0.1in
\centerline{
\hfill
{\includegraphics[width=0.47\textwidth]{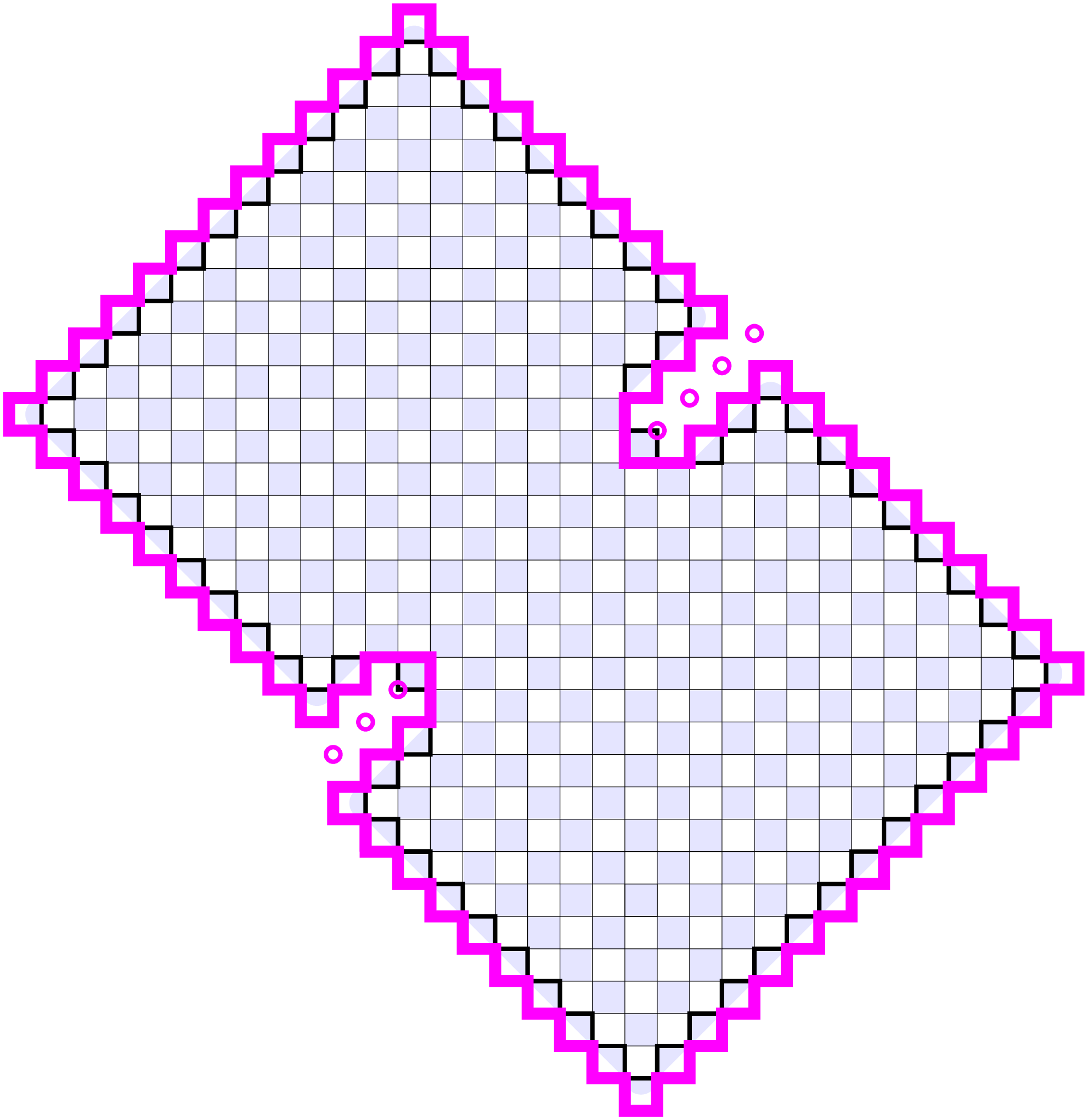}}
\hfill
{\includegraphics[width=0.47\textwidth]{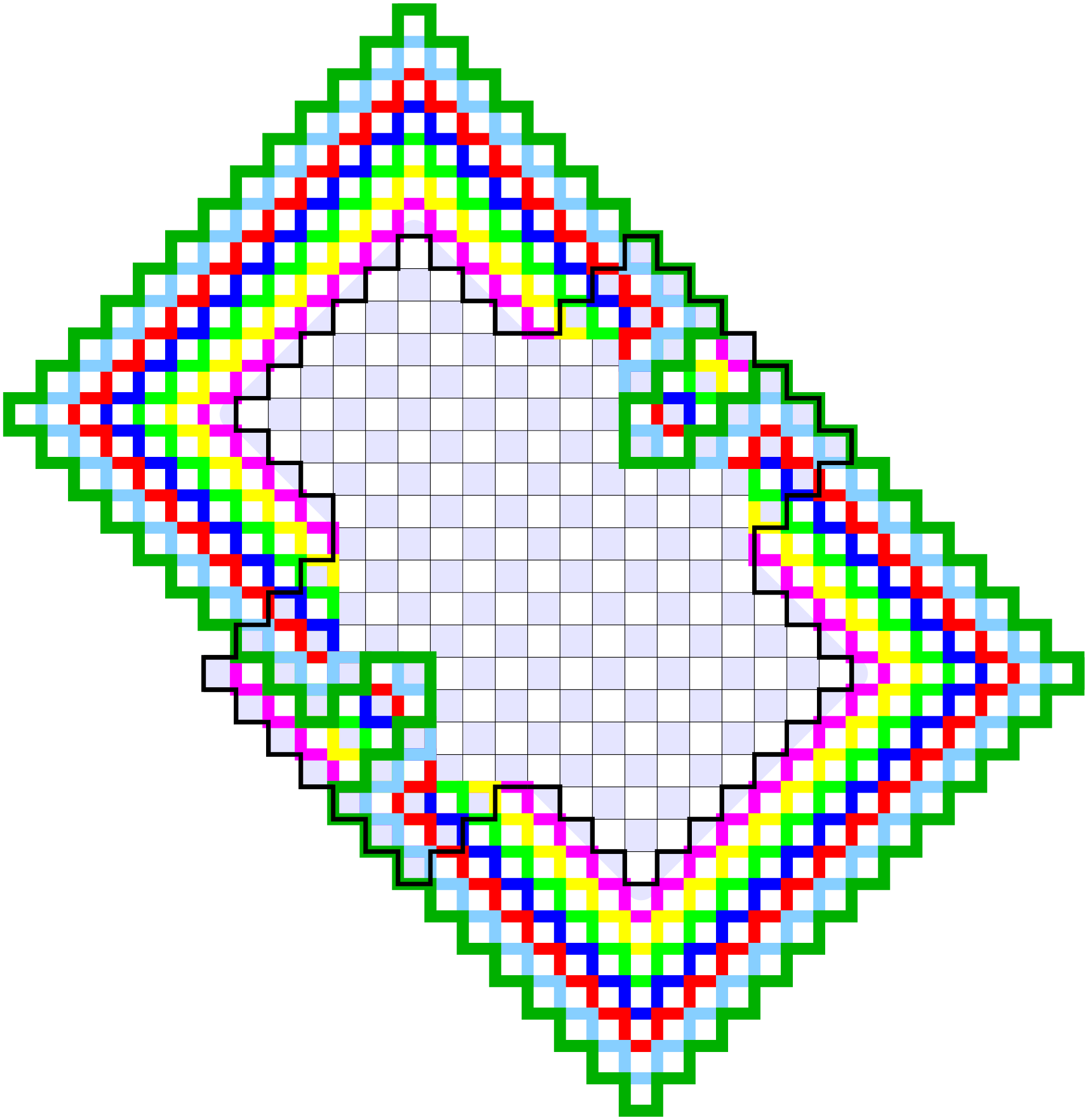}}
\hfill
}
\vskip-0.1in
\caption{
Three more applications of the complementation theorem (top left, top right and bottom left); the circles in the bottom left figure mark missing lattice points from the final graph. The evolution of the boundary of the initial nearly-cruciform graph (bounded by the black contour) to the final graph (bounded by the dark green contour); different colors indicate the intermediary stages in this evolution process.
}
\vskip-0.1in
\label{fcb}
\end{figure}

\begin{figure}[t]
\centerline{
{\includegraphics[width=0.47\textwidth]{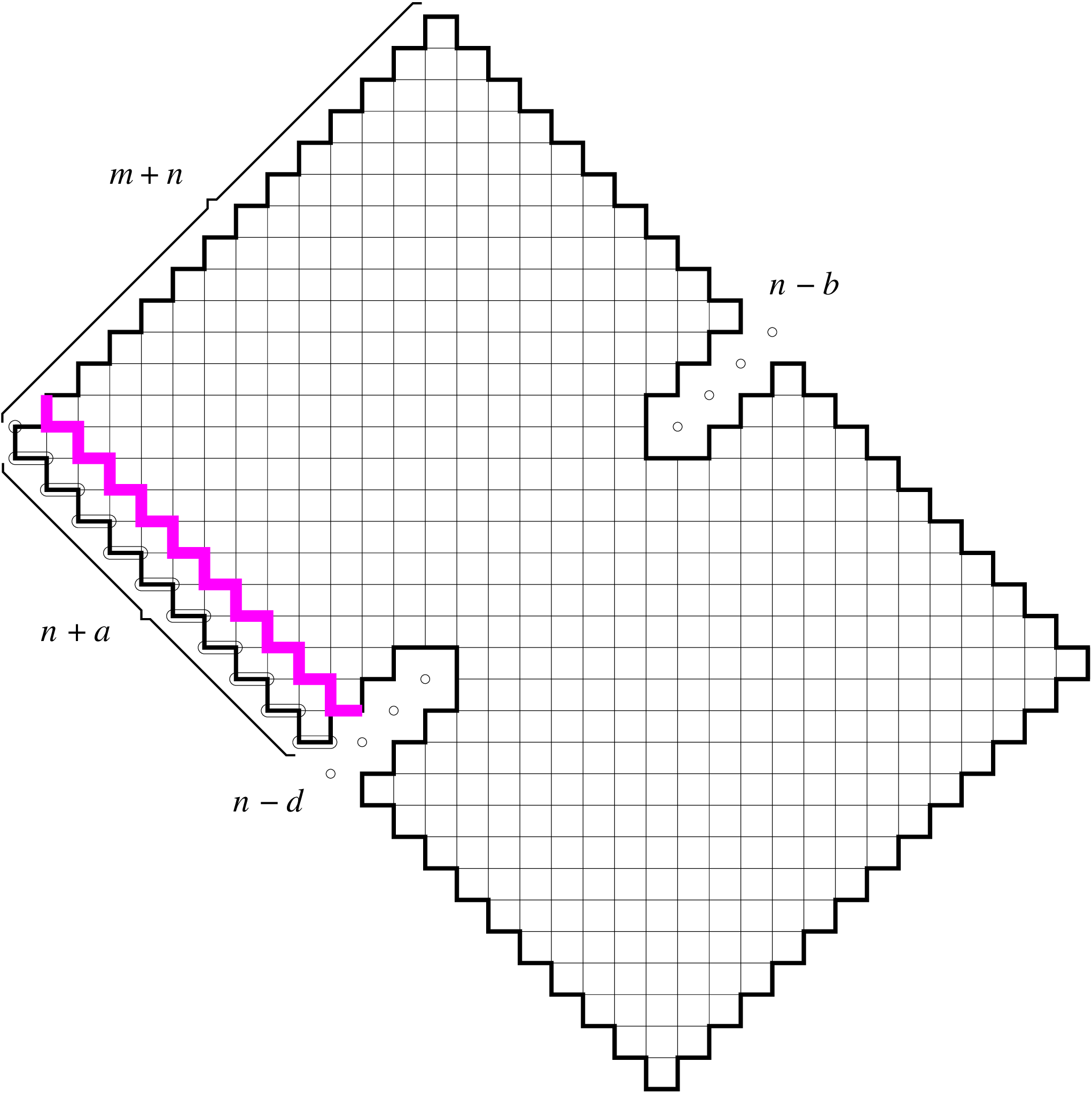}}
\hfill
{\includegraphics[width=0.47\textwidth]{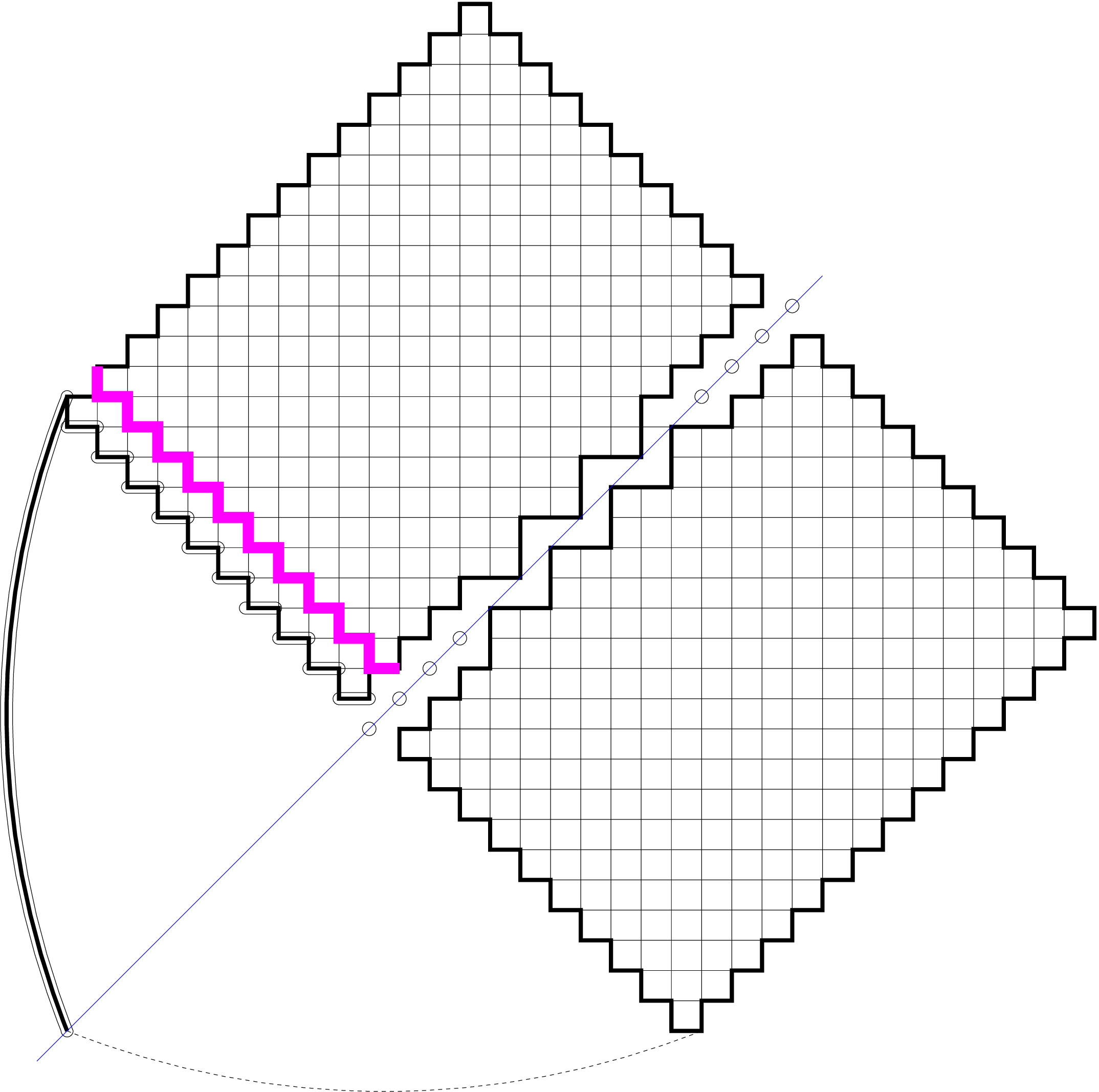}}
}
\vskip-0.1in
\caption{
{\it Left.}
$D_{m+n,0}^{a+n,b-n,a+n,d-n}$ has the same number of matchings as $\dot{AR}_{2n+2a+1,m+n}(n-d,n-b)$, the shown doubly-intruded Aztec rectangle with top left vertex removed.
{\it Right.}
When we apply the factorization theorem to a graph with twice as many matchings as $\dot{AR}_{2n+2a+1,m+n}(n-d,n-b)$, both the resulting $G^+$ and the resulting $G^-$ are $\overline{AR}_{m,n}(T)$ graphs; in this example we obtain --- after removing the forced edges --- the graph $\overline{AR}_{11,13}(1,2,3,4,6,8,10,11,12,13,14)$ (above the symmetry axis), and a graph isomorphic to $\overline{AR}_{11,14}(1,2,3,4,6,8,10,12,13,14,15)$ (below the symmetry axis).}
\vskip-0.1in
\label{fcc}
\end{figure}

What makes this proof work --- and this was the case in \cite[Theorem 2.1]{Ciu2022} as well, for the graphs considered there --- is (1) the fact that if $H$ is taken to be the nearly-cruciform graph $D_{m,n}^{a,b,a,d}$, then the complement $H'$ is precisely the graph $D_{m+1,n-1}^{a+1,b-1,a+1,d-1}$ (this is visible throughout the pictures in Figures \ref{fca} and \ref{fcb}), and (2) the graph $D_{m+n,0}^{a+n,b-n,a+n,d-n}$ has the same number of matchings as the doubly-intruded Aztec rectangle with one vertex removed described at the end of the first paragraph of this proof (this is apparent from the picture on the bottom left in Figure \ref{fcb} and the one on the left in in Figure \ref{fcc}).


The details are as follows. Apply the complementation theorem with $H$ taken to be the nearly-cruciform graph $D_{m,n}^{a,b,a,b}$. As shown on the top left in Figure \ref{fca} ($H$ is the graph bounded by the black contour and $H'$ is indicated by the thick magenta contour), the complement $H'$ is precisely $D_{m+1,n-1}^{a+1,b-1,a+1,b-1}$. A straightforward analysis of the types of the paths of cells in $H$ (see the first two pages of \cite[Section 3]{Ciu2022} and the fourth paragraph of the proof of Theorem 2.1 there) shows that the exponent $t$ of 2 in \eqref{ecb} comes out to be $n-2a-2$. Thus we obtain
\begin{equation}
\M(D_{m,n}^{a,b,a,d})=2^{n-2a-2}\M(D_{m+1,n-1}^{a+1,b-1,a+1,b-1}).
\label{ecc}
\end{equation}
Continue applying the complementation theorem successively, each time to the newly produced nearly-cruciform graph. We obtain that
\begin{equation}
\M(D_{m+i,n-i}^{a+i,b-i,a+i,d-i})=2^{n-2a-3i-2}\M(D_{m+i+1,n-i-1}^{a+i+1,b-i-1,a+i+1,d-i-1}),
\label{ecd}
\end{equation}
for $i=1,\dotsc,n-1$. Putting these together with \eqref{ecc}, we get
\begin{equation}
\M(D_{m,n}^{a,b,a,d})=2^{n(n-2a-2)-3n(n-1)/2}\M(D_{m+n,0}^{a+n,b-n,a+n,d-n}).
\label{ece}
\end{equation}
However, due to forced edges, the graph $D_{m+n,0}^{a+n,b-n,a+n,d-n}$ has the same number of matchings as the doubly-intruded Aztec rectangle with top left vertex removed $\dot{AR}_{2n+2a+1,m+n}(n-d,n-b)$ (the latter is shown on the left in Figure \ref{fcc}, where the bottom and top intrusions have lengths $n-d$ and $n-b$, respectively). Together with \eqref{ece}, this implies
\begin{equation}
\M(D_{m,n}^{a,b,a,d})=2^{n(n-2a-2)-3n(n-1)/2}\M(\dot{AR}_{2n+2a+1,m+n}(n-d,n-b)).
\label{ecf}
\end{equation}
Let $G$ be the symmetric graph obtained from $\dot{AR}_{2n+2a+1,m+n}(n-d,n-b)$ by
including back in it the deleted top left vertex and its two incident edges, as well as
including a new vertex and two new edges as shown on the right in Figure \ref{fcc}. When the factorization theorem is applied to $G$, the resulting $G^+$ and $G^-$ graphs are both graphs of the type covered by Lemma~\ref{tca}. Namely, we obtain
\begin{align}
2\M(\dot{AR}_{2n+2a+1,m+n}(n-d,n-b))
&=\M(G)
\nonumber
\\[10pt]
&=2^{m-a}\M(\overline{AR}_{n+a+1,m+n-1}(S))\M(\overline{AR}_{n+a+1,m+n}(T)),
\label{ecg}
\end{align}
where\footnote{ We use here the notation $[n]=\{1,\dotsc,n\}$.}
%
\begin{equation*}
S=[m+n]\setminus\{n-d+1,n-d+3,\ldots,n-d+2m-2a-3\}
\end{equation*}
(note that $n-d+2m-2a-3=m+b-1$, as $a+b+a+c=m+n-2$ from the balancing condition) and
\begin{equation*}
T=[m+n+1]\setminus\{n-d+1,n-d+3,\ldots,n-d+2m-2a-1\}
\end{equation*}
(by the same token, $n-d+2m-2a-1=m+b+1$).
Combining equations \eqref{ecf} and \eqref{ecg}, and using formula \eqref{eca} for the quantities on the right hand side of \eqref{ecg}, we obtain, after some manipulation, formula \eqref{ebca}.
$\hfill\square$






\end{document}